\documentclass[12pt]{amsart}

\usepackage{amsmath, amsfonts, amscd,  amssymb}
\usepackage[all]{xy}
\numberwithin{equation}{section}
\newtheorem{theorem}{Theorem}[section]
\newtheorem{proposition}[theorem]{Proposition}
\newtheorem{lemma}[theorem]{Lemma}

\theoremstyle{definition}
\newtheorem{definition}[theorem]{Definition}

\theoremstyle{remark}
\newtheorem{ass}[theorem]{Assumption}

\def \crpr {\rtimes}
\def  \nuint {\raise10pt\hbox{$\nu$}\kern-6pt\int}

\newcommand\cl{\operatorname{cl}}
\newcommand\ind{\operatorname{ind}}
\newcommand\Ind{\operatorname{Ind}}
\newcommand\ch{\operatorname{Ch}}

\renewcommand\Im{\operatorname{Im}}

\newcommand\E{\mathcal E}
\def \P{\mathcal P}

\newcommand\wm{\widehat{M}}
\newcommand\we{\widehat{E}}
\newcommand\Q{\mathcal Q}
\newcommand\R{\mathcal R}

\newcommand\B{\mathcal B}

\newcommand\D{\mathcal D}
\newcommand\Di{D\kern-6pt/}
\newcommand\cDi{{\mathcal D}\kern-6pt/}
\newcommand\CC{\mathbb C}
\def \n {\noindent}
\def \s {\smallskip}
\def \b {\bigskip}
\def \m {\medskip}
\def \cal {\mathcal}
\def \C {{\cal C}}
\def \K {{\cal K}}
\def \A {{\cal A}}
\def \Cr { C(T) \rtimes_r \Gamma}

\newcommand\NN{\mathbb N}
\newcommand\QQ{\mathbb Q}
\newcommand\RR{\mathbb R}
\newcommand\ZZ{\mathbb Z}

\newcommand\HH{\mathbb H}

\newcommand\pa{\partial}

\newcommand\Id{\operatorname{Id}}

\newcommand{\id}{\operatorname{Id}}

\newcommand{\comment}[1]                      
{
{{\bf Comment: } {\ttfamily #1}}
}

\theoremstyle{remark}
\begin{document}

\title
{Cut-and-Paste on Foliated Bundles}
\author[Eric Leichtnam]{Eric Leichtnam}
\address{Institut de Jussieu et CNRS,
Etage 7E,
175 rue du Chevaleret,
75013, Paris,  France}
        \email{leicht@math.jussieu.fr}
\author[Paolo Piazza]{Paolo Piazza}
        \address{Dipartimento di Matematica G. Castelnuovo,
Universit\`a di Roma ``La Sapienza'', P.le Aldo Moro 2, 00185 Rome, Italy}
       \email{piazza@mat.uniroma1.it}


\thanks{Both authors were partially supported by a
CNR-CNRS cooperation project. Most of this research was carried out while the
second author  was visiting  {\it Institut de Math\'ematiques de
Jussieu}, {\it Laboratoire d' Alg\`ebres d'Op\'erateurs et Repr\'esentations}. 
The financial support for this visit was provided by a {\it
Bourse de recherche de la Ville de Paris}; we heartily thank 
the {\it Ville de Paris} for having made this visit possible.}



\begin{abstract}
We discuss the behaviour
of the signature index class of closed foliated bundles
under the operation of  cutting and pasting. Along the way we establish
several index theoretic results: we define Atiyah-Patodi-Singer ($\equiv$ APS)
index classes for Dirac-type operators on foliated bundles with boundary;
we prove a relative index theorem for the difference of two APS-index classes
associated to different boundary conditions;
we prove a gluing formula
on closed foliated
bundles that are the union of two foliated bundles with boundary;
we establish a variational formula for APS-index
classes of a 1-parameter family of Dirac-type operators on foliated
bundles (this formula involves the noncommutative spectral flow of the boundary family).
All these formulas take place in the $K$-theory of a suitable cross-product
algebra.
We then apply these results in order to find sufficient conditions
ensuring the equality of the signature index classes of two cut-and-paste
equivalent  foliated bundles. We give applications to the question of when
the Baum-Connes higher signatures of closed foliated bundles are cut-and-paste invariant.

\end{abstract}

\maketitle \tableofcontents

\section{Introduction}\label{intro}

We recall that two oriented  manifolds $M_1$, $M_2$ are cut-and-paste
equivalent if
\begin{equation}\label{SK-0}
M_1=M_+ \cup_{(F,\phi_1)} M_-\,,\quad M_2=M_+ \cup_{(F,\phi_2)}
M_-
\end{equation}
with $\;\partial M_+=F=-\partial M_- \;\; \text{and}\;\; \phi_j\in
{\rm Diffeo}^+ (F).
$
In other  words, $M_1$ and $M_2$  are obtained by gluing two manifolds
with boundary but the gluing diffeomorhisms are different.

The signature of a manifold is a cut-and-paste invariant:
$\sigma(M_1)=\sigma(M_2)$ for $M_1$, $M_2$ as above. An {\it
analytic} proof of this fact is given in the book \cite{Bo-Woj}
of Booss-Bavnbek and Wojciechowski. The
argument given there is a consequence of a more general formula
concerning
 the numerical
indeces of two Dirac-type  operators obtained one from the other
by a cut-and-paste construction; the formula   expresses the
difference of the numerical indeces in terms of the spectral flow
of a suitable 1-parameter family of operators on $F$. For the
particular case of the signature operator this spectral flow turns
out to be zero, as it is simply the spectral flow of a 1-parameter
family $\{D_F(\theta)\}_{\theta\in S^1}$ of odd signature
operators on $F$ parametrized by a path of metrics. In this
vanishing result the cohomological significance of the
zero-eigenvalue has been used. Since the signature of a manifold
is equal to the index of the signature-operator, we obtain finally
\begin{equation}\label{SK-sign}
\sigma(M_2)-\sigma(M_1)={\rm ind}D_2 - {\rm ind}D_1= {\rm
sf}(\{D_F(\theta)\}_{\theta\in S^1})=0
\end{equation}
which is what we claimed.

\s
 In  this paper we shall investigate
to what extent the cut-and-paste invariance of the numerical index
of the signature operator can be generalized to the signature
index class  of  {\it foliated bundles}, i.e. $\Gamma$-equivariant
fibrations $\widehat{M}\to T$ with $T$ a manifold on which
$\Gamma$ acts, and $\widehat{M}$ a manifold on which $\Gamma$-acts
freely properly and cocompactly; thus the quotient
$\widehat{M}/\Gamma :=M$ is  a smooth manifold and
$\Gamma\rightarrow
 \widehat{M}\rightarrow M$
is a $\Gamma$-Galois covering.

We shall be therefore interested in index theory on foliated
bundles, both in the closed case and in the case where a boundary
is present. The relevant index classes, for Dirac type operators,
will live in $K_0 (C(T)\rtimes_r \Gamma)$ if the fibers of the
$\Gamma$-equivariant fibration $\widehat{M}\to T$ are
even-dimensional and in $K_1 (C(T)\rtimes_r \Gamma)$ if the fibers
are odd-dimensional; $C(T)\rtimes \Gamma$ denotes the reduced
cross-product algebra. Notice that the manifold $\widehat{M}/\Gamma$
is foliated by the images of the fibers of the fibration under the
projection map $\widehat{M}\rightarrow \widehat{M}/\Gamma$. These
foliations can be quite interesting; in fact it is well known that
one can get any type (I, II, III) of foliation for suitable
choices of $\Gamma$-equivariant fibrations. Index theory on
foliated bundles is a particular but important case of the general
foliation-index-theory developed by  Connes, see \cite{Co}, and Connes-Skandalis
\cite{CoSka}.

Notice that if $T=$ point and $\Gamma=\{1\}$ then we simply have a
compact manifold $M$. Moreover $C(T)\rtimes_r \Gamma=\CC$ and (in
the even dimensional case) the index class is nothing but the
numeric index of the operator under the isomorphism $K_0
(\CC)=\ZZ$.
 If $\Gamma=\{1\}$ we simply have a fibration, $C(T)\rtimes_r \Gamma = C(T)$ and the
index class
 reduces to the Atiyah-Singer family  index in $K_* (C(T))=K^* (T)$.
 Finally, if $T=$ point then we have a Galois
covering, $C(T)\rtimes_r \Gamma=C^*_r \Gamma$, the reduced group
$C^*$-algebra associated to $\Gamma$, and
 the index class, which now lives in $K_* (C^*_r \Gamma)$,
 is nothing but the Mishchenko-Fomenko index class associated
to the Dirac operator twisted by the canonical flat line bundle of
the covering.

 In previous work of
ours, with collaborators, we investigated the cut-and-paste
invariance of the signature index class  in the case of Galois
$\Gamma$-coverings $\Gamma\to \widetilde{M}\to M$, thus solving
(at least partially) a problem raised by Lott and also by
Weinberger \cite{Lott 3}. See Leichtnam-Lott-Piazza \cite{LLP} for
the first positive results in this direction and then
Leichtnam-Lueck-Kreck \cite{LLK} and Leichtnam-Piazza
\cite{LPCUT}. It was explained in \cite{LPCUT} that for Galois
$\Gamma$-coverings the signature index class, in $K_* (C^*_r
\Gamma)$, is not cut-and-paste invariant: one shows that the
difference of signature index classes for two cut-and-paste
equivalent coverings
\begin{equation}\label{intro:difference} \Ind (\D^{{\rm
sign}}_{\widetilde{M}_2})-\Ind (\D^{{\rm sign}}_{\widetilde{M}_1})
\end{equation}
is equal to a {\it higher} spectral flow, in $K_* (C^*_r \Gamma)$,
for an $S^1$-family of $C^*_r \Gamma$-linear signature operators
on the cutting hypersurface $F$. This formula is the consequence
of \begin{itemize}
\item a gluing formula for the index class of a closed Galois covering
which is the union of two Galois coverings with boundary;
\item a variational formula for the index classes associated
to a path of Dirac operators on a Galois covering with boundary
\end{itemize}

The definition of higher spectral flow was given by Dai and Zhang
\cite{Dai-Zhang} for a path of families of Dirac operators
parametrized by a compact speace $T$, i.e. for a path of
$C(T)$-linear operators; this definiton is based on the notion of
spectral section, given by Melrose and Piazza in \cite{MP I}. The
papers \cite{LPGAFA} and \cite{LPCUT} extend the results of
Melrose-Piazza and Dai-Zhang from the family-case, i.e.
$C(T)$-linear operators, to the Galois-coverings case, i.e. $C^*_r
(\Gamma)$-linear operators. This step should be thought of as the
passage from a commutative to a non-commutative context.

In contrast with the numeric case explained above, the higher
spectral flow appearing in formula (\ref{intro:difference}) will
{\it not} be equal to zero, in general. It is however possible to
give sufficient conditions on the cutting hypersurface $F$
ensuring the vanishing of this higher spectral flow and therefore
the equality of the two signature index classes. This hypothesis
comes from Lott's paper \cite{Lott II}. If, in addition, the group
$\Gamma$ satisfies the Strong Novikov Conjecture (i.e. the
rational injectivity of the assembly map), then the equality of
the index classes implies the equality of the Novikov higher
signatures \footnote{We recall that for a Galois covering
$\Gamma\to \widetilde{M}\to M$ the Novikov higher signatures are
the numbers: $$\int_M L(M) \cup r^* [c]\,,\quad [c]\in H^*
(B\Gamma,\CC)=H^* (\Gamma,\CC)$$ with $r:M\to B\Gamma$ the
classifying map of the covering.}. These ideas are now explained
in the survey of Leichtnam-Piazza \cite{LPBOUTET}.

\n
Summarizing:

\begin{itemize}
\item suitable conditions on the cutting
hypersurface $F$ ensures that the signature index class on Galois
$\Gamma$-coverings is a cut-and-paste invariant;
\item further conditions on the group $\Gamma$ allow
to deduce the cut-and-paste invariance of all Novikov higher
signatures from the cut-and-paste invariance of the signature
index classes
\end{itemize}

In the present paper we wish to follow the above line of reasoning
for the more general case of foliated bundles. The specific
problems we wish to solve are the following:
\begin{itemize}
\item  give sufficient conditions ensuring that the signature index
class, in the group $K_* (C(T)\rtimes_r \Gamma)$, is a cut-and-paste
invariant;
\item find additional conditions on $\Gamma$ and its action on $T$
ensuring that the Baum-Connes higher signatures (a generalization
to foliated bundles of the Novikov higher signatures) are
cut-and-paste invariant.
\end{itemize}
In order to solve the first problem we will need to develop a
general Atiyah-Patodi-Singer index theory on foliated bundles.
 Some of our
arguments will be easy extensions of the Galois coverings case and  we
will be quite brief in such cases; other arguments will be more
involved and we shall explain them in detail. We shall use the
$\Gamma$-equivariant $b$-pseudodifferential calculus on foliated
bundles with boundary developed in \cite{LPETALE} where a
Atiyah-Patodi-Singer index theory was developed under an
invertibility assumption on the boundary operator.

\bigskip
{\bf The paper is structured as follows.} In {\bf Section \ref{sect:sflow} } we
recall the arguments leading to the cut-and-paste invariance of
the {\it numeric} index of the signature operator: conceptually
this is the model case that will be extended to our more general
situation. In {\bf Section \ref{sect:foliated-bundles}} we begin by
recalling the definition of index class associated to a
$\Gamma$-equivariant family of Dirac-type operators
$(D(\theta))_{\theta\in T}$ on a foliated bundle $\wm\to T$; we
denote by $\D$ the $C(T)\rtimes_r \Gamma$-linear operator defined
by the family $ (D(\theta))_{\theta\in T}$. We then introduce the
notion of spectral section associated to $\D$ and prove the
fundamental existence theorem: {\it a spectral section $\P$ for
$\D$ exists if and only if the index class $\Ind(\D)$ in $K_*
(C(T)\rtimes_r \Gamma)$ vanishes.} We also introduce the notion of
difference class $[\P]-[\Q]$, in $K_{*+1} (C(T)\rtimes_r \Gamma)$,
associated to two spectral sections $\P, \Q$; following Dai and
Zhang we then introduce the notion of higher spectral flow for a
path $(\D_u)_{u\in [0,1]}$ of $\Cr$-linear operators. In {\bf Section
\ref{sect:b-aps-classes}} we develop index theory on foliated
bundles {\it with boundary}, using the $b$-pseudodifferential
calculus on foliated bundles developed by Leichtnam and Piazza in
\cite{LPETALE}; thus we start in Subsection
\ref{subsect:preliminaries} by  reviewing the numeric case,
explaining the equality between the {\it generalized} APS-index on
a manifold with boundary  and a certain {\it perturbed}
$L^2$-index on the associated manifold with cyclindrical ends. The
latter can also be described in the framework of Melrose'
$b$-geometry \cite{Melrose}; this will be in fact the point of
view that we shall adopt. In Subsection
\ref{subsect:foliated-with-boundary} we describe in detail the
geometric  set-up for foliated bundles with boundary; we also
recall the index theory developed  in \cite{LPETALE} for
$\Gamma$-equivariant families of Dirac-type operators
$(D(\theta))_{\theta\in T}$ with invertible boundary family. We
explain how spectral sections for the boundary family can be used
in order to remove the invertibility assumption; we then prove in
Subsection \ref{subsect:cobordism} the cobordism invariance of the
index class in this general case; by the existence theorem we
infer that a boundary family always admits a spectral section. In
Subsections \ref{subsect:b} and \ref{subsect:APS} we define $b$ and
APS-index classes in $K_* (C(T)\rtimes_r \Gamma)$ associated to a
$\Gamma$-equivariant family $(D(\theta))_{\theta\in T}$ {\it and}
a choice of spectral section for the boundary family. We also
prove the equality of these two index classes.  In {\bf Section
\ref{sect:fundamentals}} we establish 3 fundamental properties of
these index classes: {\it the gluing formula} for index classes on
closed foliated bundles that are union of two foliated bundles
with boundary; {\it the relative index theorem}, equating the
difference of  $b$-index classes associated to two different
choices of spectral sections $\P$, $\Q$, to the  class
$[\Q]-[\P]$; {\it the variational formula} computing the variation
of the $b$-index of a  path of $\Gamma$-equivariant families in
terms of the higher spectral flow associated to the boundary
family. Section \ref{sect:b-aps-classes} and
\ref{sect:fundamentals} are modeled on the work of Melrose-Piazza
\cite{MP I} \cite{MP II} and the subsequent work of
Leichtnam-Piazza \cite{LPGAFA} \cite{LPCUT}. In {\bf Section
\ref{sect:cut-and-paste}} we finally tackle the cut-and-paste
problem: we define two cut-and-paste equivalent foliated bundles
in Subsection \ref{subsect:Cut-and-paste}. We then compute the
difference of index classes for the signature family of two
cut-and-paste equivalent manifold in terms of the higher spectral
flow of a path of operators on the cutting hypersurface
(Subsection \ref{subsect:defect}); we refer to such a formula as a
{\it defect formula}. In Subsection \ref{subsect:vanishing-higher}
we employ a gap condition on forms of middle degree
(see Assumption \ref{Lottgroupoid}) and spectral
sections with a certain symmetry property in order to give
conditions ensuring the vanishing of this defect; this will be a
solution to the problem we had posed. Finally, in {\bf Section
\ref{sect:applications}} we give additional conditions on $\Gamma$
and its action on $T$ in order to deduce the cut-and-paste
invariance of certain geometric numerical invariants generalizing
the Novikov higher signatures. One of our main geometric results in Section
 \ref{sect:applications} is the following:

{\it Assume that the rational Baum-Connes map
$$
 \mu_\QQ: K_{0,\tau} ((E\Gamma\times T)/\Gamma)\otimes_\ZZ \QQ \rightarrow K_0 ( \Cr)\otimes_\ZZ \QQ.
 $$
is injective.
Let $\widehat{X}_\phi\to T$ and $\widehat{X}_\psi\to T$ be
two $\Gamma$-equivariant fibrations that are cut-and-paste equivalent 
and satisfy  
 Assumption \ref{Lottgroupoid} below.
Assume  moreover that the  vertical tangent bundles both admit a $\Gamma$-invariant spin structure.
Then for any
$c \in  H^{\ast}( (E\Gamma \times T)/\Gamma ; \QQ)$ the Baum-Connes higher signatures
are equal:
$$\int_{\widehat{X}_\phi
/{\Gamma} }
L( \widehat{X}_\phi
/{\Gamma})\wedge r_\phi^*(c)\,=\,\int_{\widehat{X}_\psi
/\Gamma }
L(\widehat{X}_\psi
/\Gamma )\wedge s_\psi^*(c).
$$
}

{\bf Acknowledgements.} 
We thank Alexander Gorokhovsky, John Lott, George Skandalis and Jean Louis Tu 
for helpful discussions.

\section{Spectral flow}\label{sect:sflow}

\subsection{Spectral flow through spectral sections.}
\label{subsect:sflow+ss}

Let $(D_t)_{t\in [0,1]}$ be a continuous family of formally
self-adjoint elliptic differential operators. For simplicity, we
shall
 assume that  $D_0$ and $D_1$ are
invertible. The spectral flow of the family $(D_t)_{t\in [0,1]}$,
$${\rm sf}((D_t)_{t\in [0,1]})\in \ZZ\,,$$
 is by
definition the net number of eigenvalues changing sign as $t$ runs
from $0$ to $1$.

Following Dai and Zhang \cite{Dai-Zhang} we shall now recall how
to express the spectral flow of the family $(D_t)_{t\in [0,1]}$ in
a way which can be generalized to situations where the spectrum is
not discrete. To this end we recall the notion of spectral section
associated to the family $ (D_t)_{t\in [0,1]}$ and, in fact, to
any family of formally self-adjoint elliptic differential
operators parametrized by a compact space $B$. Thus, let
$\mathcal{D}=(D_z)_{z\in B}$ be a smooth family of formally
self-adjoint elliptic differential operators parametrized by $B$.
 We shall eventually take $B=[0,1]$ but let us
proceed in full generality for the time being.
 Each operator
$D_z$ acts on the sections of a hermitian vector bundle $F_z$ over
a closed riemannian manifold $N_z$.  A spectral section $\P$ for
$\D$ is a {\bf smooth} family $\P=(P_z)_{z\in B}$ of self-adjoint
projections with $P_z\in\Psi^0(N_z;F_z)$ and satisfying the
following property:
\begin{equation}\label{ss}
\exists\, R\in\RR\; \;|\;\;D_z u=\lambda u
\Rightarrow Pu=u\;\;
\text{if}\;\;\lambda>R,\,\quad Pu=0\;\;\text{if}\, \lambda<-R
\end{equation}
This means that each operator $P_z$ is equal to the identity on
the eigenfunctions of $D_z$ corresponding to large positive
eigenvalues and equal to 0 on those corresponding to large
negative eigenvalues. Each $P_z$ is a finite rank perturbation of
the spectral projection $$\Pi_{\geq} (z):= \chi_{[0,\infty)}
(D_z)$$ corresponding to the non-negative eigenvalues of $D_z$.
The family $\D$ is a family of self-adjoint Fredholm operators
parametrized by $B$; it therefore defines an index class
$\Ind(\D)\in K^1(B)$. One of the main results of Melrose-Piazza
\cite{MP I}   asserts
that a spectral section for $\D$ exists if and only if
$\Ind(\D)=0$ in $K^1(B)$. One can prove that if $\Ind(\D)=0$ then
the set of spectral sections associated to $\D$ is infinite.

Let us go back to our 1-parameter family of formally self-adjoint
elliptic differential operators
 $\D=(D_t)_{t\in [0,1]}$ on an odd dimensional manifold $N$;
 since $B=[0,1]$ is contractible, it is
 certainly the case that
$\Ind(\D)=0$. Thus there exists a  spectral section
$\P=(P_t)_{t\in [0,1]}$ for the 1-parameter family $\D=(D_t)_{t\in
[0,1]}$. Consider now the spectral projection $\Pi_{\geq} (0)$
associated to the nonnegative eigevalues  of $D_0$; consider the
spectral projection $\Pi_{\geq} (1)$ associated to the nonnegative
eigenvalues of $D_1$; consider the projection $P_0$ and the
projection $P_1$, also associated, respectively, to $D_0$ and
$D_1$. Since $P_0$ is a finite rank perturbation of $\Pi_{\geq}
(0)$ one can define a relative index $i(\Pi_{\geq} (0),P_0)\in
\ZZ$; this is simply the index of the Fredholm operator
$\Pi_{\geq} (0)\circ P_0: \Im P_0\to \Im \Pi_{\geq} (0)$ .
Similarly, there is a relative index $i(\Pi_{\geq} (1) ,P_1)\in
\ZZ$; we shall denote these relative indeces as $[\Pi_{\geq} (0)
-P_0]$, $[\Pi_{\geq} (1) -P_1]$ respectively

\begin{proposition} (Dai-Zhang \cite{Dai-Zhang})\\
If $\D=(D_t)_{t\in [0,1]}$ is a smooth 1-parameter family of
formally self-adjoint elliptic differential operators, with $D_0$
and $D_1$ invertible, then
\begin{equation}
{\rm sf}((D_t)_{t\in [0,1]})=[\Pi_{\geq} (1)  -P_1] -[\Pi_{\geq}
(0) -P_0]\in \ZZ\,.
\end{equation}
\end{proposition}
We can see the spectral flow as an element in $K^0 ({\rm
point})=K_0 (\CC)=\ZZ$ (the K-theory of an algebra will be defined
in Subsection  \ref{subsect:fredholm} below).

This result leads in a natural way to a small generalization of
the notion of spectral flow: fix  spectral sections $Q_i$ for
$D_i$, $i=0,1$. If $\P=(P_t)_{t\in [0,1]}$ is a total spectral
section for the family $(D_t)_{t\in[0,1]}$, then the spectral flow
${\rm sf}((D_t)_{t\in [0,1]}; Q_0,Q_1)$ from $(D_0,Q_0)$ to
$(D_1,Q_1)$ through $(D_t)_{t\in [0,1]}$ is the element of $K_0
(\CC)$ given by the difference class
\begin{equation}\label{sf-general}
{\rm sf}((D_t)_{t\in [0,1]}; Q_0,Q_1)\,:=\,[Q_1 -
P_1]-[Q_0-P_0]\;\;\in\;\;K_0 (\CC)=\ZZ\;;
\end{equation}
one can prove that this class is well defined, independent of the
total spectral section $\P=(P_t)_{t\in [0,1]}$ chosen. The classic
case explained above is obtained by making the particular choice
$Q_0= \Pi_{\geq} (0)$, $Q_1=\Pi_{\geq} (1)$.

\subsection{Index and spectral flow}\label{subsect:index-and-ss}

Let $M$ be an even-dimensional riemannian manifold with boundary,
endowed with a product metric near the boundary . In contrast with
the closed case, the Atiyah-Patodi-Singer ($\equiv$ APS) index of
a Dirac-type operator, see \cite{APS}, is not stable under
perturbations.
 In fact, assume  that $(D_t)_{t\in [0,1]}$
is a smoothly varying family of Dirac operators on $M$; as an
important example we could consider a family of metrics
$(g(t))_{t\in [0,1]}$ on $M$ and the associated family of
signature operators $(D^{{\rm sign}} (t))_{t\in [0,1]}$. We could 
also consider, on a spin manifold, a family of 
Dirac operators $(\Di (t))_{t\in [0,1]}$ parametrized by a path 
of metrics $(g(t))_{t\in [0,1]}.$
Going
back to the general case, consider the family of operators induced
on the boundary $(D_{\pa M} (t))_{t\in [0,1]}\,$; let $\Pi_{\geq}
(t)$ the corresponding spectral projection associated to the
non-negative eigenvalues. For simplicity, let us assume that the
boundary operator is invertible at $t=1$ and at $t=0$; then the
following variational formula for the APS-indeces holds:
\begin{equation}\label{variational}
\ind (D^+ _1,\Pi_{\geq} (1))-\ind (D^+ _0,\Pi_{\geq} (0))={\rm sf}
((D_{\pa M} (t))_{t\in [0,1]})\,.
\end{equation}
Formula (\ref{variational}) follows from the APS-index formula.
It can also be proved analytically, without
making use of the APS-index formula. See for example
\cite{Dai-Zhang} where much more general projections are allowed.
In fact, in that paper Dai and Zhang establish a
more general variational formula; since such a generalization will
be important to us we briefly explain it.

First of all, if $D$ is an odd Dirac-type operator on $M$, acting
on the sections of a $\ZZ_2$-graded Clifford module $E=E^+\oplus E^-$,
 and if $Q$
is a spectral section for $D_{\pa M}$, then there is a well
defined {\it generalized} APS-boundary value problem, with a
well-defined index $\ind
(D^+,Q)$ (see for example \cite{Bo-Woj}); the boundary problem
is simply defined by taking the operator $D^+$ with domain
$$\{u\in C^\infty (M,E^+)\;|\; u|_{\pa M}\in {\rm Ker} Q\}\,.$$
This generalized boundary
value problems has interesting properties. First of all, if
$Q^\prime$ is a different spectral section for $D_{\pa M}$, then
the following {\it relative index formula} holds:
\begin{equation}\label{relative-baby}
\ind (D^+,Q^\prime)-\ind (D^+,Q)=[Q-Q^\prime]\in K_0 (\CC)=\ZZ\,.
\end{equation}
Second, let $(D_t)_{t\in [0,1]}$ be a smoothly varying family of
odd Dirac operators on $M$; choose a spectral section $Q_0$ for
$D_0$ and a spectral section $Q_1$ for $D_1$; then the following
{\it variational formula for generalized APS-indeces} holds:
\begin{equation}\label{variational-general}
\ind (D^+ _1,Q_1)-\ind (D^+ _0,Q_0)={\rm sf} ((D_{\pa M}
(t))_{t\in [0,1]};Q_1,Q_0)\,.
\end{equation}

\subsubsection{{\bf Remark.}}\label{subsubsect:important}
 If $N$ is odd dimensional and
$(D^{{\rm sign}}_N (t))_{t\in [0,1]}$ is a one-pa\-ra\-me\-ter family of
odd signature operators parametrized by a path of metrics $g_N
(t)_{t\in [0,1]}$, then
\begin{equation}\label{zero-spectral-flow}
{\rm sf}((D^{{\rm sign}}_N (t))_{t\in
[0,1]};\Pi_{\geq}(0),\Pi_{\geq}(1))=0\,,
\end{equation}
\begin{equation}
{\rm sf}((D^{{\rm
sign}}_N (t))_{t\in [0,1]};\Pi_{>}(0),\Pi_{>}(1))=0 \,.
\end{equation}
In fact, the kernel of the odd signature operator is equal to the
space of  harmonic forms on $N$; from the Hodge theorem we know
that such a vector space is independent of the metric we choose;
thus there are not eigenvalues changing sign and we can choose as
total spectral section $\P= (\Pi_{\geq}(t))_{t\in [0,1]}$ for the
first equation and $\P= (\Pi_{>}(t))_{t\in [0,1]}$ for the second
equation. In particular, if $(D^{{\rm sign}}_M (t))_{t\in [0,1]}$
is a one-parameter family of signature operators (parametrized by
a path of metrics) on a 4k-dimensional manifold with boundary,
then we have
\begin{equation}\label{constant-index}
\ind (D^{{\rm sign},+}_M (1),\Pi_{\geq} (1))=\ind (D^{{\rm
sign},+}_M (0),\Pi_{\geq} (0))\end{equation}
\begin{equation}\label{constant-index-bis}
 \ind (D^{{\rm sign},+}_M (1),\Pi_{>}
(1))=\ind (D^{{\rm sign},+}_M (0),\Pi_{>} (0))\,.
\end{equation}

\subsection{The gluing formula.}\label{subsect:gluing-numeric}

We consider $X$, a {\it closed} oriented
compact manifold which is the union of two manifolds with
boundary. Thus there exists an embedded hypersurface $F$ which
separates $M$ into two connected components and such that
$$X=M_+\cup_F M_-\,,\quad\text{with}\quad \pa M_+=F=-\pa M_-\,.$$
We assume that the metric $g$ is of product type near the
hypersurface $F$, i.e. near the boundaries of $M_+$ and $M_-$. Let
$D_X$ be a Dirac-type operator on $X$; then we obtain in a natural
way two Dirac operators on $M_+$ and $M_-$. The following gluing
formula holds:
\begin{equation}\label{additivity-simple}
\ind (D_X)=\ind (D_{M_+},\Pi_{\geq})+ \ind (D_{M_-},1-\Pi_{\geq})
\,. \end{equation} The discrepancy in the spectral projections
comes from the orientation of the normals to the two boundaries (if
one is inward pointing, then the other is outward pointing).
\footnote{Notice that $ 1-\Pi_{\geq}$ is not exactly the
APS-projection associated to the non-negative eigenvalues of
$D_{\pa M_-}$; to be precise  $ 1-\Pi_{\geq}=\Pi^{\pa M_-}_{>},$
the projection onto the {\it positive} eigenvalues of $D_{\pa
M_-}$\/.} In particular, for the signature operators we have the
fundamental formula
\begin{equation}\label{additivity-simple-bis}
\ind (D^{{\rm sign}}_X)=\ind (D^{{\rm sign}}_{M_+},\Pi_{\geq})+
\ind (D^{{\rm sign}}_{M_-},1-\Pi_{\geq}) \,. \end{equation}

\smallskip

The two formulae (\ref{additivity-simple})
(\ref{additivity-simple-bis}) can be proved directly, in a purely
analytical fashion, see Bunke \cite{Bunke}, Leichtnam-Piazza \cite{LPCUT}. Of course they
are  also a consequence of the APS-index theorem.

\subsection{An analytic proof of the cut and paste invariance of the signature.}

The gluing formula (\ref{additivity-simple-bis}) for the signature
operator
 can be
generalized to a more complicated situation, where $X_\phi$ is a
closed manifold obtained by {\it gluing}  two manifolds with
boundary  through a diffeomorphism $\phi$ between their
boundary. We shall concentrate on the signature operator. Thus let
$M$ and $N$ be two oriented manifolds with boundary and let
$\phi:\pa M\to \pa N$ be an oriented diffeomorphism. We consider the manifold
with boundary $X_\phi:= M\cup_\phi N^-$, with $N^-$ equal to $N$
with the opposite orientation.
We shall follow the notation of the previous
subsection;  thus we set $M_+ := M$, $M_- := N^-$ and 
$X_\phi = M_+ \cup_\phi M_- .$

We fix a metric $g_\phi$ on $X_\phi$.
 Notice that giving $g_\phi$ is equivalent to give
$g(+)$ on $M_+$ and $g(-)$ on $M_-$ such that $\phi^*(g(-)|_{\pa
M_-})=g(+)|_{\pa M_+}.$ We shall assume that these metrics are of
product type near the boundary. The pull-back $\phi^*$ defines an
isometry between $L^2(\pa M_-,\Lambda^* (\pa
M_-))$ defined by $ g(-)|_{\pa M_-}$ and $L^2 (\pa
M_+,\Lambda^* (\pa M_+))$ defined by $g(+)|_{\pa M_+}$. Let
$D^{{\rm sign}}_{X_\phi}$ be the signature operator on $X_\phi$
associated to $g_\phi$. We also have the signature operators on
$M_\pm$ with boundary operators $D^{{\rm sign}}_{\pa M_+}$ and
$D^{{\rm sign}}_{\pa M_-}$ and it is easy to check that
\begin{equation}\label{conj}
D^{{\rm sign}}_{\pa M_+}= - \phi^* (  D^{{\rm sign}}_{\pa M_-} )
(\phi^*)^{-1} \,.
\end{equation}

 Let  $\Pi^{\pa M_+}_{\geq}$ be  the APS  spectral
projection for $ D^{{\rm sign}}_{\pa M_+}$ and consider the
projection $$\Pi^\phi_{\geq}:=(\phi^*)^{-1}  \Pi_{\geq} \phi^*
\,;$$ from (\ref{conj}) we infer  that  ${\rm Id}-
\Pi^\phi_{\geq}$ is equal to the spectral projection $\Pi_{>}^{\pa
M_-}$ onto the non-negative eigenvalues of $D^{{\rm sign}}_{\pa
M_-}$. Here the fact that we are dealing with the signature
operator has been used.  One
can prove, analytically, the following additivity formula:
\begin{align*}
\ind (D^{{\rm sign}}_{X_\phi}) &= \ind (D^{{\rm
sign}}_{M_+},\Pi^{\pa M_+}_{\geq}) + \ind (D^{{\rm
sign}}_{M_-},1-\Pi^\phi_{\geq})\\& =\ind (D^{{\rm
sign}}_{M_+},\Pi^{\pa M_+}_{\geq}) + \ind (D^{{\rm
sign}}_{M_-},\Pi_{>}^{\pa M_-} )\,.
\end{align*}

Summarizing, also in this more general case we have
\begin{equation}\label{additivity-bis}
\ind (D^{{\rm sign}}_{X_\phi})=\ind (D^{{\rm sign}}_{M_+},\Pi^{\pa
M_+}_{\geq}) + \ind (D^{{\rm sign}}_{M_-},\Pi_{>}^{\pa M_-} )\,,
 \end{equation}
where it is important to notice that the operators appearing in
this formula are associated to the metrics $g_\phi$ on the left
hand side and to  metrics $g_\phi (+)$, $g_\psi (-)$ on the right
hand side; thus we should write more precisely
\begin{equation}\label{additivity-bis-precise}
\ind (D^{{\rm sign}}_{(X_\phi,g_\psi)})=\ind (D^{{\rm
sign}}_{(M_+,g_\psi (+))},\Pi^{\pa M_+}_{\geq}) + \ind (D^{{\rm
sign}}_{(M_-,g_\psi (-))},\Pi_{>}^{\pa M_-} )\,.
 \end{equation}

 Let now $$X_\phi=M_+ \cup_\phi M_- \,,\quad
X_\psi=M_+\cup_\psi M_-\,$$ with $\phi,\, \psi: \pa M_+
\longrightarrow \pa M_-$ diffeomorphisms, be two such manifolds.
One says in this case that $X_\phi$ and $X_\psi$ are {\it
cut-and-paste equivalent}. Let us fix metrics $g_\phi$ and
$g_\psi$ on $X_\phi$ and $X_\psi$ respectively. We obtain metrics
$g_\phi(\pm)$, $g_\psi (\pm)$ on $M_\pm$. We can assume these
metrics to be product like near $\pa M_\pm$. Let $D^{{\rm
sign}}_{X_\phi}$ and $D^{{\rm sign}}_{X_\psi}$ be the signature
operators associated to $g_\phi$ and $g_\psi$. We wish to give a
proof of the equality $$\ind (D^{{\rm sign},+}_{X_\phi}) = \ind
(D^{{\rm sign},+}_{X_\psi})\,.$$ It will be important to keep
track of the metrics involved, thus we write as above $D^{{\rm
sign}}_{(X_\phi,g_\phi)}$ for the signature operator on $X_\phi$
associated to the metric $g_\phi$ and \linebreak$D^{{\rm
sign}}_{(M_\pm,g_\phi (\pm))}$ for the induced signature operators
on the manifold with boundary $M_\pm$. Similarly we proceed for
 $$D^{{\rm sign}}_{(X_\psi,g_\psi)} \quad\text{and}\quad
 D^{{\rm sign}}_{(M_\pm,g_\psi(\pm))}$$
 We
begin by applying the additivity formula: we obtain
\begin{align*}
 \ind (D^{{\rm sign}}_{(X_\phi,g_\phi)})=&
\ind (D^{{\rm sign}}_{(M_+,g_\phi (+))},\Pi^{\pa M_+}_{\geq}) +
\ind (D^{{\rm sign}}_{M_-,g_\phi (-)}, \Pi^{\pa M_-}_{>})\\ \ind
(D^{{\rm sign}}_{(X_\psi,g_\psi)})=& \ind (D^{{\rm
sign}}_{(M_+,g_\psi(+)) },\Pi^{\pa M_+}_{\geq}) + \ind (D^{{\rm
sign}}_{(M_-,g_\psi (-) )}, \Pi^{\pa M_-}_{>})\,.
\end{align*}
On the left hand side of these formulae we have indeces of
operators on manifolds which are, in general, non-diffeomorphic.
On the right hand side, on the other hand, we can compare, as we
have the same 2 manifolds with boundary. From the above formula we
infer that
\begin{align*}
& \,\ind (D^{{\rm sign},+}_{(X_\psi,g_\psi)})- \ind (D^{{\rm
sign},+}_{(X_\phi,g_\phi)})\\ &= \left(\ind (D^{{\rm
sign}}_{(M_+,g_\psi (+)) },\Pi^{\pa M_+}_{\geq})- \ind (D^{{\rm
sign}}_{(M_+,g_\phi (+))},\Pi^{\pa M_+}_{\geq})\right)\\ &\, +
\left(\ind (D^{{\rm sign}}_{(M_-,g_\psi (-)) },\Pi^{\pa M_-}_{>})-
\ind (D^{{\rm sign}}_{(M_-,g_\phi (-))},\Pi^{\pa M_-}_{>})\right)\,.
\end{align*}
Let $(g_t (\pm))_{t\in [0,1]}$ be a path of metrics on $M_\pm$
joining $g_\phi (\pm)$ and $g_\psi (\pm)$. Applying the
variational formula (\ref{variational}) to the two summands on the
right hand side  we obtain:
\begin{align}\label{variational-cut}
& \,\ind (D^{{\rm sign},+}_{(X_\psi,g_\psi)})- \ind (D^{{\rm
sign},+}_{(X_\phi,g_\phi)})\\ &= {\rm sf} \left( (D^{{\rm
sign}}_{(\pa M_+, g_t (+))})_{t\in [0,1]}; \Pi^{\pa M_+}_{\geq}
(0),\Pi^{\pa M_+}_{\geq} (1) \right)\\ &\, + {\rm sf} \left(
(D^{{\rm sign}}_{(\pa M_-, g_t (-))})_{t\in [0,1]}; \Pi^{\pa
M_-}_{>} (0), \,,\Pi^{\pa M_-}_{>} (1) \right)
\end{align}
and taking into account Remark \ref{subsubsect:important}, we
immediately obtain $$\ind (D^{{\rm sign},+}_{(X_\psi,g_\psi)})=
\ind (D^{{\rm sign},+}_{(X_\phi,g_\phi)})\,;$$
this implies the
equality of the signatures as required. In fact a small argument
involving our various identifications shows that our conclusion
can be reached through the following two equalities
\begin{equation}\label{zero-with-mapping-torus}
\Ind (D^{{\rm sign},+}_{(X_\psi,g_\psi)})- \Ind (D^{{\rm
sign}}_{(X_\phi,g_\phi)})= {\rm sf} (\{D^{\rm sign}_{{\rm odd}}
(\theta)\}_{\theta\in S^1})=0 \,.\end{equation} The spectral flow
appearing in this formula is associated to a $S^1$-family of odd
signature operators acting on the fibers of the mapping torus $M(F
, \phi^{-1} \circ \psi)\to S^1$ and parametrized by a family of
metrics. As already remarked 
this
spectral flow is zero because of the cohomological significance of
the zero eigenvalue for the signature operator.

\s
\n
{\bf Remark.} It should be remarked that in this analytic proof of
the cut-and-paste invariance of the signature, we have not used
the APS-index formula; only the analytic properties of the APS
boundary value problem were employed. This will be important
later, when we shall extend the argument above to the foliated
case.

\section{Foliated bundles, index classes and
the noncommutative spectral flow}\label{sect:foliated-bundles}

\subsection{Preliminaries: K-Theory of $C^*$-algebras and Fredholm operators.}
\label{subsect:fredholm}

Let $A$ be a unital $C^*$-algebra. We
recall that $K_0 (A)$ is defined as the Grothendieck group
associated to the semigroup  of isomorphism classes of
 finitely generated projective left
 $A-$modules. $K_0 (A)$
is an additive group.
A class in $K_0 (A)$ is represented by a formal difference
$[E]-[F]$ of isomorphism classes of finitely generated projective
left $A$-module. Notice that a finitely generated projective
left $A$-module is the range of a projection $p$ in the matrix algebra
$M_n(A)$, for a suitable $n$. In fact, $K_0 (A)$ can be also described
in terms of such projections: one considers the inductive limit
$M_\infty (A)$, the cartesian product $M_\infty (A)^2$
and identifies two pairs of projections $(p,q) \in
M_n(A)^2$ and $(p^\prime, q^\prime)\in M_{n^\prime} (A)^2$ if for
suitable $k, k^\prime \in \NN$, $$p\oplus  q^\prime \oplus {\rm
Id}_k \oplus 0_{k^\prime}\;\;\text{is conjugate to}\;\; p^\prime
\oplus q \oplus {\rm Id}_k \oplus 0_{k^\prime}\;\;\text{in}\;\;
M_{n+n^\prime+k+k^\prime}(A ).$$ The quotient under this
equivalence relation has a natural structure
of abelian group and is  naturally isomorphic to $K_0 (A)$;
one denotes by $[p]-[q]$
(=$[p^\prime]-[q^\prime]$) the class of $(p,q)$.
Recall that if $(p_1,q_1)\in M_{n_1} (A)^2$ then one has:
$([p]-[q]) + ([p_1]-[q_1])=[p\oplus p_1] -[q\oplus q_1]$
where $([p\oplus p_1], [q\oplus q_1]) \in M_{n+n_1}(A)^2.$
 When $A$ is a non unital $C^*-$algebra one
introduces the unital $C^*-$algebra $\widetilde{A}= A \oplus \CC$
obtained by adding the unit element $0 \oplus 1$ to $A$; one
considers the morphism $\epsilon: \widetilde{A} \rightarrow \CC$
defined by $\epsilon (a \oplus \lambda) = \lambda.$ One then
defines $K_0(A)$ to be equal to the kernel of the map $\epsilon_*:
K_0(\widetilde{A}) \rightarrow K_0(\CC)$ induced by $\epsilon$.
Observe that $K_0(\CC)= K_0( M_n(\CC) ) = \ZZ.$ We  define
$K_1(A)$ to be equal to $K_0( A \otimes C_0(\RR) )$ where $ A
\otimes C_0(\RR)$ is the {\it suspension of $A$}. For instance
$K_1(\CC)= K_1(M_n(\CC) )= 0.$ Alternatively,  $K_1(A)$ can be
identified with the set of connected components of $GL_\infty(A).$
Fundamental properties of the $K_*$-functor are the Bott
isomorphism $$K_0 (A)\simeq K_0 (C_0 (\RR^2)\otimes A)$$ and
the six-terms long
exact sequence associated to the short exact sequence
$$0\longrightarrow J\overset{i} \longrightarrow A \overset{\pi}
\longrightarrow A/J
\longrightarrow 0$$
with $J$ an ideal in $A$:


\def\mapver#1{\Big\downarrow\rlap{$\vcenter{\hbox{$\scriptstyle#1$}}$}}
\def\vermap#1{\Big\uparrow\rlap{$\vcenter{\hbox{$\scriptstyle#1$}}$}}

\[
\begin{matrix}
& K_0 (J) & \overset{i_*}\longrightarrow & K_0 (A)& \overset{\pi_*}\longrightarrow & K_0(A/J) &\\
& \vermap{\partial} & & & &  \mapver{\partial} &\\
& K_1(A/J)& \overset{\pi_*}\longleftarrow & K_1 (A) & \overset{i_*}\longleftarrow & K_1 (J)&
\end{matrix}
\]

\m
Elements in $K_* (A)$ arise naturally as index classes
of generalized Fredholm operators between Hilbert $A$-modules.
We recall that a Hilbert $A$-module $\mathcal{E}$ is a left
$A$-module endowed with an $A$-valued form $\langle\,,\,\rangle:\mathcal{E}
\times\mathcal{E}\to A$ satisfying the following axioms:

(i) $\langle \eta,\xi_1+\xi_2\rangle=\langle\eta,\xi_1 \rangle +\langle\eta,\xi_2\rangle$;

(ii) $\langle\eta,\xi a\rangle=\langle\eta,\xi\rangle a$;

(iii) $\langle\eta,\xi\rangle ^*=\langle\xi,\eta\rangle$;

(iv) $\langle\eta,\eta\rangle \geq 0$;

(v) $\langle\eta,\eta\rangle=0$ $\Leftrightarrow$ $\eta=0$

(vi) $\mathcal{E}$ is complete with respect to the norm $\| \eta \|:=
\|\langle\eta,\eta\rangle\|_A$

\m
Hilbert $A$-modules share {\it some} of the usual properties
of Hilbert spaces; there are however structural differences
and more care is needed in all arguments involving
bounded operators, adjoints, orthogonal complements etc...
See for example \cite{WO} for a clear account of
these structural differences. Here we shall content ourselves
by simply stating
that there is indeed a natural notion of bounded
$A$-linear adjointable operator between
Hilbert $A$-modules \footnote{
Notice that a bounded $A$-linear operator might not admit an adjoint.}.
One defines in this way the space of bounded adjointable
operators between two Hilbert $A$-modules
$\mathcal{E}^\pm$, denoted $\mathcal{B}_A(\mathcal{E}^+,
\mathcal{E}^-)$. There is also a natural
notion of finite rank operator and the closed subspace
of compact operators,
$\mathcal{K}_A (\mathcal{E}^+,
\mathcal{E}^-)\subset \mathcal{B}_A(\mathcal{E}^+,
\mathcal{E}^-)$, is defined as the norm closure in
$\mathcal{B}_A(\mathcal{E}^+,
\mathcal{E}^-)$
of the space of finite rank operators.

 A $A$-Fredholm operator $\mathcal{L}^+: \mathcal{E}^+
\longrightarrow \mathcal{E}^-$ is a bounded adjointable
operator which
is invertible modulo compacts.  One can prove that {\it up to a compact
perturbation} the kernel of $\mathcal{L}^+$ and of its adjoint
are finitely generated projective $A$-modules; we get in this way
an index class $\Ind (\mathcal{L}^+)\in K_0 (A)$ by simply
taking the formal difference
of these finitely generated projective $A$-modules.
If $\mathcal{L}:\mathcal{E}\to \mathcal{E}$ is self-adjoint, then there
is an index class in $K_1 (A)$. Let us describe it. First we need
to recall another fundamental property of the $K_*$-functor, namely
its {\it stability}: {\it if $\mathcal{E}$ is a Hilbert $A$-module,
then $K_* (\mathcal{K}_A (\mathcal{E}))\simeq K_* (A)$}.
Consider now the Calkin algebra $$\mathcal{C}_A (\mathcal{E}):=
\frac{\mathcal{B}_A(\mathcal{E})}{\mathcal{K}_A (\mathcal{E})}\;$$
and the canonical projection map $\pi: \mathcal{B}_A(\mathcal{E})
\to \mathcal{C}_A (\mathcal{E})$;
the long exact sequence in K-theory associated to
$$ 0\rightarrow
\K_{A} (\mathcal{E})\rightarrow \B_{A}(\mathcal{E})
\rightarrow \C_A (\mathcal{E})\rightarrow 0 $$
gives the homorphism  $$ \delta:
K_0(\C_{A}(\mathcal{E})) \rightarrow
K_1(A)\,,$$
since, by stability,
$K_1(\K_{A}(\mathcal{E}))\simeq K_1(A).$
We need to  assume that $ \mathcal{L}^2 =\id + \mathcal{R}$ where
$\mathcal{R}\in \K_{A}$; then
  $\pi (\frac{1}{2}(\mathcal{L}+\Id))$ is a projection
 in the Calkin algebra and
 the index class $\Ind (\mathcal{L})$ is defined by
 $$\Ind (\mathcal{L}):=\delta [\pi (\frac{1}{2}(\mathcal{L}+\Id))]
\in K_1(A). $$

\subsection{Foliated bundles.}\label{subsect:foliated-bundles}

Let $\Gamma$ be a finitely generated discrete group. Let $T$ be a
smooth closed compact connected manifold on which $\Gamma$ acts on
the right. Let $\widehat{M}$ be a closed manifold on which
$\Gamma$ acts freely, properly and cocompactly on the right: the
quotient space $M=\widehat{M} /\Gamma$ is thus a smooth closed
compact manifold. We assume that $\widehat{M}$ fibers over $T$ and
that the resulting fibration $$ \pi: \,\widehat{M} \rightarrow T
$$ is a {\it $\Gamma$-equivariant fibration} with  fibers
$\pi^{-1}(\theta),\theta \in T,$.

\s
\n
{\bf Remark.} We observe incidentally that what we have described
is an example of a proper cocompact $G$-manifold $P$ with $G$ an
\'etale groupoid, see Connes ( \cite{Co} page 137) for the definition. In
our case $$G=T\rtimes\Gamma\,,\quad\quad G^{(0)}=T\,,\quad\quad
(\alpha:P\to G^{(0)})\equiv (\pi:\wm \to T).$$ The  groupoid $G =
T \crpr \Gamma$ has  $G^{(1)}=T\times \Gamma$ as set of morphisms
and $G^{(0)}= T$ as base. The range and source maps are
respectively given by: $$ \forall (\theta,g) \in T \times
\Gamma,\; r(\theta,g)=\theta,\;\;\; s(\theta,g)= \theta\cdot g\,.
$$ The composition is defined as follows: $$(\theta,g)\cdot
(\theta^\prime,g^\prime)=(\theta,g
g^\prime)\;\;\text{if}\;\;\theta^\prime=\theta g\,.$$ The inverse
of $(\theta,g)$ is  $(\theta g,g^{-1})$.

\s
With a small abuse of terminology
we shall call  $\Gamma$-equivariant fibrations $\widehat{M}\to T$
with smooth compact quotient $M=\widehat{M}/\Gamma$
a  proper $T\rtimes \Gamma$-manifold. It is important to notice that
the compact manifold $M$ inherits a foliation ${\cal
F}$, with leaves equal to the images of the fibers of $\pi:
\,\widehat{M} \rightarrow T $ under the quotient map
$\widehat{M}\rightarrow M=\widehat{M}/\Gamma$. The foliated
manifold $(M,{\cal F})$ is usually referred to as a {\it foliated
$T$-bundle} or simply as a {\it foliated bundle}.

\s
{\it In this paper we shall refer to a $\Gamma$-equivariant
fibrations $\widehat{M}\to T$
with smooth compact quotient $M=\widehat{M}/\Gamma$ either as
a proper $T\rtimes \Gamma$-manifold or as a foliated bundle.}

\medskip
\noindent {\bf Example.} Let $X$ be a compact closed manifold  and
let $\Gamma\rightarrow \widetilde{X}\rightarrow X$ be a Galois
cover of $X$. Let $T$ be a smooth compact manifold on which
$\Gamma$ acts by diffeomorphisms. We consider
$\widehat{M}=\widetilde{X}\times T$, $ \pi=\text{projection onto
the second factor}$, $M=\widetilde{X}\times_{\Gamma}
T:=(\widetilde{X}\times T)/\Gamma$ where we let $\Gamma$ act on
$\widetilde{X}\times T$ diagonally. The leaves of the foliation
${\cal F}$ are the images of the manifolds $\widetilde{X}\times
\{\theta\}$, $\theta\in T$.

\s
\n
As a {\it particular example} of this construction consider $T=S^1$,
$\Gamma=\ZZ$, $\widetilde{X}=\RR$, so that $\widehat{M}=\RR\times S^1$.
We let $n\in \ZZ$ act on $(r,e^{i\theta})\in \RR\times S^1$ by
$n\cdot (r,e^{i\theta}):= (r+n,e^{i(\theta+n\alpha)})$, for some fixed
$\alpha\in\RR$. Then $M=T^2$
and if $\alpha/2\pi $ is irrational we get as a foliation
of $T^2$ the well-known Kronecker
foliation. This is a type II foliation.

\s
\n
As a  different example, consider a
smooth closed riemann surface $\Sigma$ of genus $g>1$ and let
$\Gamma=\pi_1 (\Sigma)$, a discrete subgroup of $PSL(2,\RR)$. Then
we can consider $X:=\Sigma$, $\Gamma\rightarrow\widetilde{X}
\rightarrow X$ equal to the universal cover
$\Gamma\to\HH^2\rightarrow \Sigma$, $T=S^1$, with $\Gamma \unlhd
PSL(2,\RR)$ acting on $S^1=\RR P^1$ by fractional linear
transformations. The resulting  foliation is of type III.

\b
Let us go back to the general situation; we shall also denote the
typical fiber of $\pi: \,\widehat{M} \rightarrow T $ by $Z$. We
now assume that $Z$ is of dimension $2k-1$. We choose a
$\Gamma$-invariant metric
 on the vertical tangent
bundle $ TZ $. Finally, we  assume the existence of a
$\Gamma-$equivariant spin structure on $ TZ $ that is fixed once
and for all. We denote  by $S^Z \rightarrow \widehat{M} $ the
associated spinor bundle. We consider also a $\Gamma-$equivariant
complex hermitian vector bundle $ \widehat{V} \rightarrow
\widehat{M}$ endowed with a $\Gamma-$invariant hermitian
connection. We then set $\widehat{E}=S^Z \otimes \widehat{V}$
which defines a smooth $\Gamma-$invariant family of hermitian
Clifford modules on the fibers $\pi^{-1}(\theta), \theta \in T.$
We thus get a $\Gamma-$equivariant family $(D(\theta))_{\theta \in
T}$ of  Dirac type operators acting fiberwise on
$C^\infty_c(\widehat{M} , \widehat{E}).$
If $Z$ is even-dimensional, then the spinor bundle $S^Z$ is $\ZZ_2$-graded;
thus $\widehat{E}$ is also $\ZZ_2$-graded, $\widehat{E}=
\widehat{E}^+\oplus\widehat{E}^-$, and the family $(D(\theta))_{\theta\in T}$
is now odd with respect to this grading:
 $$  D(\theta)= \begin{pmatrix}
0 & {D^-(\theta)} \cr {D^+(\theta)}  & 0 \cr
\end{pmatrix},\; \theta\in T.
$$
Notice that we could assume that $\we$ is defined more generally
by a smooth $\Gamma-$in\-va\-ri\-ant family of $\ZZ_2-$graded hermitian
Clifford modules on the fibers $\pi^{-1}(\theta),\, \theta \in T$.

\smallskip
{\it We shall now explain how such a $\Gamma$-equivariant
family defines a class in the K-theory
of the cross-product
algebra $C(T)\rtimes_r \Gamma$; thus we shall define
suitable Hilbert $C(T)\rtimes_r \Gamma$-modules and see
how the family $(D(\theta))_{\theta\in T}$ defines a
$C(T)\rtimes_r \Gamma$-Fredholm operator $\D$ on them.}

\subsection{The $C^*$-algebra $C(T)\rtimes_r \Gamma$}

The {\it algebraic} cross-product $C^\infty_c(T) \crpr \Gamma$ is, by
definition, the set of functions $\sum_{g \in \Gamma} t_g(\theta)
g$ such that only a finite number of the $ t_g(\cdot) \in
C^\infty_c(T)$ do not vanish identically. We shall identify any
function  $f$ having compact support $$ f: T \crpr \Gamma
\rightarrow \CC $$ $$ (\theta,g) \rightarrow f(\theta,g) $$ with $
\sum_{g \in \Gamma} f(\theta,g) g$. Then one has: $$
\sum_{g^\prime  \in \Gamma} f^\prime (\theta, g^\prime) g^\prime
\cdot \sum_{g  \in \Gamma} f (\theta, g) g\,=\, \sum_{h \in
\Gamma} \left( \sum_{g\in\Gamma} f^\prime
(\theta,g^\prime)\,f(\theta\cdot g^\prime,(g^\prime)^{-1} h)
\right) h $$  where we recall that $$g^\prime \cdot
(f(\theta,g) g) = f(\theta\cdot g^\prime,g) g^\prime g\,.$$ The
algebraic cross-product $C^\infty_c(T) \crpr \Gamma$ will also be
denoted by $C^\infty_c(T \crpr \Gamma)$. One can introduce the
reduced $C^*$-algebra $C^*_r (T\rtimes \Gamma)$ associated to the
groupoid $T\rtimes \Gamma$ as a suitable completion of the
algebraic cross product $C^\infty_c (T) \crpr \Gamma$. See
\cite{Co}.

It is well known, and easy to check, that  there is a
natural isomorphism between the reduced $C^\ast$-algebra
$C_r^\ast(T\rtimes \Gamma)$ of the groupoid $T\rtimes \Gamma$ and
the cross-product algebra $C(T)\rtimes_r \Gamma$ (see, for
example, Moore-Schochet \cite{MoS}); we shall henceforth identify these two
$C^*$-algebras. Observe, $T$ being compact,
 that the reduced
$C^*-$algebra $C(T)\rtimes_r \Gamma$ is {\it unital}.

\subsection{$C(T)\rtimes_r \Gamma$-Hilbert modules.}
\label{subsect:hilbert-modules}

Recall that the action of $\Gamma$ on $\widehat{E}$ induces for
each $g\in \Gamma$ an operator $R^*_{g}$ acting on
${C}^\infty_c(\widehat{M}, \widehat{E} )$.
We endow $C^\infty_c(\widehat{M}, \widehat{E} )$ with the structure of
 left $C^\infty_c(T) \crpr \Gamma-$module by setting
for any $s \in C^\infty_c(\widehat{M}, \widehat{E} )$ and
$\sum_{g \in \Gamma} f(. ,g) g \in C^\infty_c(T) \crpr \Gamma$
\begin{equation} \label{U}
\forall p \in \widehat{M},\; \;
\left( \sum_{g \in \Gamma} f(. ,g) g \,\cdot s \right)\, (p):=
 \sum_{g \in \Gamma} f(\pi(p), g) (R^\ast_g  s\,)(p)\,.
\end{equation}

We  define a $C^\infty_c(T) \crpr \Gamma$-valued
hermitian product
of two sections $s$ and $s^\prime$ of
${C}^\infty_c(\widehat{M}, \widehat{E} )$ by setting:
\begin{equation} \label{hermitianproduct}
\langle s;s^\prime\rangle \,=\,
\sum_{g\in \Gamma} \langle s;s^\prime\rangle (\theta,g) g\;\;\in
C^\infty_c(T) \crpr \Gamma
\subset\Cr
\end{equation}
where $\forall (\theta,g) \in T \times \Gamma $:
$$
 \langle s;s^\prime\rangle(\theta,g)\,=\,
\int_{\pi^{-1}(\theta.g)}\, \langle R^\ast_{g^{-1}} (s)(y) ; s^\prime(y)
\rangle_{\widehat{E}} d {\rm Vol}_{\pi^{-1}(\theta\cdot g)}(y) $$ with
$ d{\rm Vol}_{\pi^{-1}(\theta\cdot g)}(y)$ denoting the riemannian
density in the fiber.

\s
{\it Summarizing}: we have defined on
${C}^\infty_c(\widehat{M}, \widehat{E} )$
a structure of $C^\infty_c(T) \crpr \Gamma$-module; moreover we have
defined
a form
$$\langle\,\,,\,\,\rangle: {C}^\infty_c(\widehat{M}, \widehat{E} )
\times {C}^\infty_c(\widehat{M}, \widehat{E} )
\rightarrow C^\infty_c(T) \crpr \Gamma\subset\Cr\,.$$

One denotes by $ L^2_{C(T)\rtimes_r
\Gamma}(\widehat{M} , \widehat{E})$  the completion of
${C}^\infty_c(\widehat{M}, \widehat{E} )$ with respect to
the norm induced by
$\langle\,;\,\rangle$.
One can prove that $ L^2_{C(T)\rtimes_r
\Gamma}(\widehat{M} , \widehat{E})$ is indeed a Hilbert
$C(T)\rtimes_r
\Gamma$-module. In a similar way, one can introduce Sobolev-modules
$H^m_{C(T)\rtimes_r \Gamma}(\widehat{M} , \widehat{E}), m\in\NN$.

\subsection{Index classes on closed foliated bundles}
\label{index-classes-closed}

Let us go back to our $\Gamma$-equivariant family  of Dirac-type
operators $(D(\theta))_{\theta\in T}$.

\begin{lemma}  \label{D} The family of Dirac operators
$(D(\theta))_{\theta \in T}$ acting fiberwise
on \linebreak$C^\infty_c(\widehat{M}, \widehat{E} )$
 defines a left $C^\infty_c(T \crpr \Gamma)-$linear endomorphism $\D$ of
$C^\infty_c(\widehat{M}, \widehat{E} )$.
\end{lemma}

\begin{proof}
 Using the above notations we have:
$$ \D \bigl( \, \sum_{g \in \Gamma} f(. ,g) g \,\cdot s \,
\bigr)(p) =
 \sum_{g \in \Gamma} f(\pi(p), g) \, D(\pi(p)) (R^\ast_g  s\,) \, (p)
$$ where we have used the fact that $(D(\theta))_{\theta \in T}$
is a
 family of operators, i.e.
commutes with the natural action of $C^\infty(T)$. Since  the
family $(D(\theta))_{\theta \in T}$ is $\Gamma-$equivariant,
 the right hand  is by definition equal to
$$ \left( \sum_{g \in \Gamma} f(. ,g) g \,\cdot \D( s)\right) \,
(p) $$ which proves the lemma.
\end{proof}

One can prove that $\D$ extends to a bounded operator
$$
\D: H^m_{C(T)\rtimes_r \Gamma}(\widehat{M} , \widehat{E}) \rightarrow
 H^{m-1}_{C(T)\rtimes_r \Gamma}(\widehat{M} , \widehat{E})
 $$
for each $m\in \NN$. Moreover, using an appropriate
pseudodifferential  calculus,
  one can show that
each extension
is $\Cr$-Fredholm. See Subsection \ref{subsect:perturbations}
below for more on this point.
If the fibers are even-dimensional, then the family
$(D(\theta))_{\theta\in T}$ is $\ZZ_2$-graded odd and we get
in this way an index class $\Ind(\D^+)\in K_0 (\Cr)$
which is  independent of $m$ by elliptic regularity.
Alternatively, let $\mathcal{E}:=L^2_{C(T)\rtimes_r \Gamma}(\widehat{M} ,
\widehat{E})$;
we can consider the bounded operator
$$\mathcal{L}=\frac{\D}{(\Id+\D^2)^{\frac{1}{2}}}\;;$$
the fact that $\mathcal{L}$
is well defined in $\B_{C(T)\rtimes_r \Gamma}(\mathcal{E})$,
the algebra of
bounded operators on the
Hilbert $\Cr$-module $\mathcal{E}$, requires proof and is based on the
continuous functional calculus for {\it regular} unbounded
operators on Hilbert modules. The reader interested in  the details 
may read Proposition 1 of \cite{LPETALE} and Baaj-Julg's work 
explained in \cite{Co} page 433.
If the fibers are even dimensional, then $\Ind(\D^+)\equiv \Ind (\mathcal{L}^+)$.

If the fibers are odd dimensional we define the index  class
in $K_1(C(T)\rtimes_r \Gamma)$ proceeding as in Subsection
\ref{subsect:fredholm}.
 Thus, by definition,
$$\Ind (\D)=\delta [\pi(\frac{1}{2}(\mathcal{L}+\Id))]\in K_1(\K_{C(T)\rtimes_r
\Gamma}(\mathcal{E}))\simeq K_1(C(T)\rtimes_r \Gamma)$$
with $\K_{C(T)\rtimes_r \Gamma} (\mathcal{E})\subset
\B_{C(T)\rtimes_r \Gamma}(\mathcal{E})$ denoting the sub-algebra of
$C(T)\rtimes_r \Gamma$-compact operators and
$$\pi: \B_{C(T)\rtimes_r \Gamma}(\mathcal{E})
\rightarrow \C_{C(T)\rtimes_r \Gamma}(\mathcal{E})=
\frac{\B_{C(T)\rtimes_r \Gamma}(\mathcal{E})}
{\K_{C(T)\rtimes_r \Gamma} (\mathcal{E})}$$
denoting the projection onto the Calkin algebra.

 \subsection{Spectral sections for Dirac operators on  foliated bundles.}

Both the notion of spectral flow and the
APS-boundary value problem are based on the possibility
of ``dividing
in two parts'' the spectrum of a self-adjoint Dirac operator
on a closed manifold.
This is done via the self-adjoint projection
$\chi_{[0,\infty)}(D)$
associated to the non-negative
spectrum of the relevant operator $D$. In the noncommutative context
things are more complicated. Thus let $\D$ be a $C(T)\rtimes_r \Gamma$-linear
Dirac operator, associated to a $\Gamma$-equivariant family.
As already remarked there is a well defined {\it continuous}
functional calculus associated to $\D$; however, as we are working
in a $C^*$-algebraic context, there is {\it not} a measurable
calculus; thus it does not make sense to consider the operator
$\Pi_{\geq} (\D):=\chi_{[0,\infty)} (\D)$ as an element in
$\B_{C(T)\rtimes_r \Gamma}(L^2_{C(T)\rtimes_r \Gamma}(\widehat{M} ,
\widehat{E}))$.\footnote{ Needless to say, if $\D$ is invertible, from
$H^1_{\Cr}(\wm ,\we)$ onto  $L^2_{\Cr}(\wm ,\we)$,
then (see Proposition 1 of
\cite{LPETALE})
we can indeed define $\Pi_{\geq} (\D)$, either by taking a smooth approximation
of the characteristic function $\chi_{[0,\infty)}$ or by setting
$$\Pi_{\geq} (\D) := \frac{1}{2} \left( \frac{\D}{|\D|}+\Id
\right)\,.$$} This is the reason why we need
a more general notion of spectral projection
for ``dividing  in two parts'' the spectrum of
$\D$; this  is provided by the definition of {\it  spectral section}.
We have already encountered this notion for families of Dirac operator
parametrized by a compact manifold $T$, see \ref{subsect:sflow+ss};
such a family defines a $C(T)$-linear operator and our task now
is to pass from the commutative case to the noncommutative case
where a group $\Gamma$ is present and the relevant linearity
is with respect to the noncommutative $C^*$-algebra $C(T)\rtimes_r \Gamma$.

\s
In what follows we shall briefly denote the relevant algebras
of $C(T)\rtimes_r \Gamma$-linear operators by
$$\B_{C(T)\rtimes_r \Gamma}\,,\quad \mathcal{K}_{C(T)\rtimes_r \Gamma}\,,
\quad \C_{C(T)\rtimes_r \Gamma}\,.$$

\begin{definition} A spectral cut $\chi$ is by definition a function
$\chi \in C^\infty( \RR, [0, 1])$ such that $\chi(x)=0$ for $x
<<0$ and $\chi(x)=1$ for $x >>1.$
\end{definition}

Observe that $\chi(\D)$ induces
a projection in the Calkin algebra
$\C_{C(T)\rtimes_r \Gamma}$
which does not depend on the choice of the spectral cut $\chi.$
In fact, always in $\C_{C(T)\rtimes_r \Gamma}$, we have the equality
$$\pi(\frac{1}{2}(\mathcal{L}+\Id))=\pi(\chi(\D))\,,\text{with}\;\;
\mathcal{L}=\frac{\D}{(\Id+\D^2)^{\frac{1}{2}}}\;.$$
Thus the index class $\Ind \D \in K_1(C(T)\rtimes_r \Gamma ) $
 is also defined by $\Ind \D= \delta [\chi (\D)]$ for
 any spectral cut $\chi$.

 \begin{definition} A spectral section
 $\P$ for $\D$ is a self-adjoint projection $\P \in \B_{C(T)\rtimes_r \Gamma}$
such  that there exist
two spectral cuts $\chi_1, \chi_2$ such that $\chi_2 \equiv 1$
on a neigborhood of the support of $\chi_1$ and
 ${\rm Im}\,\chi_1(\D) \subset  {\rm Im}\,\P \subset {\rm Im}\,\chi_2(\D).$
\end{definition}

\begin{theorem}\label{theo:existence} {\item 1)} If $\D$ admits a spectral section then
$\Ind \D=0$ in $K_1(C(T)\rtimes_r \Gamma ) .$ {\item 2)} Assume
that $\Ind \D=0$ in $K_1(C(T)\rtimes_r \Gamma ) .$ Then $\D$
admits a spectral section $\P$.
 \end{theorem}
\begin{proof}

1) The proof of \cite{LPCUT} page 363 extends immediately to this
context.

2) The proof of Theorem 3 of \cite{LPCUT} (partially based
on the unpublished work of Wu \cite{Wu}) shows that
we just have to prove that $\id- \chi(\D)$ and $\chi(D)$
define, for any spectral cut $\chi$, two very full projections of
$\C_{C(T)\rtimes_r \Gamma}.$ Then the proof
of Lemma 4 of \cite{LPCUT} (page 360) shows that we just have to prove
that for a given spectral cut $\chi_1$ there exists
$u \in L^2_{C(T)\rtimes_r \Gamma}(\widehat{M} , \widehat{E})$
such that $ \langle \chi_1(\D)(u) ; \chi_1(\D)(u) \rangle$ is
invertible in $C(T)\rtimes_r \Gamma$.

There exists
a positive integer $N$ and open connected subsets
$U_i \subset U_i^\prime$ of $T$ ($1 \leq i \leq N$) with the following
properties. Each $U_i$ is relatively compact in $U_i^\prime$,
$\cup_{1\leq i \leq N} U_i = T $, and for each $i\in \{1,\ldots,N\}$
the restriction of the fibration $\pi$ to
$\pi^{-1}(U_i^\prime)$ is trivial: $\pi^{-1}(U_i^\prime) \simeq
U_i^\prime \times Z$. Denote by $\mu$ a given $\Gamma-$invariant
 riemannian measure on $\widehat{M}$ and $v_f$ the volume of a fundamental domain
 for the action of $\Gamma$ on $\wm$. Then using an induction
 argument on $i\in \{1, \ldots,N\}$, we may find an open connected
 subset $W \subset Z$ and  open
 subsets $V_i^\prime \subset \pi^{-1}(U_i^\prime)$  ($1 \leq i \leq N$)
 with the following three properties:

 \n (a)  $ V_i^\prime \simeq U_i^\prime \times W$ for each $i\in \{1, \ldots,N\}$
  and $ \mu( \cup_{1 \leq i \leq N} V_i^\prime ) \leq
 {1 \over 2} v_f$

 \n (b) $ \forall \gamma \in \Gamma \setminus \{e\}, \,
 \forall i \in \{1,\ldots , N \},\; \overline{V_i^\prime}
 \cdot \gamma \cap \overline{V_i^\prime} =
 \emptyset$.

\n (c) $ \forall \gamma \in \Gamma, \,
 \forall i,j \in \{1,\ldots , N \},\,
 {\rm { with}}\, i \not=j,\; \overline{V_i^\prime}\cdot
 \gamma \cap \overline{V_j^\prime} =
 \emptyset$.
 Now, if the $U^\prime_i$ are small enough then the proof of
 Lemma 4 of \cite{LPCUT} shows that for any $\epsilon>0$ and
 $i \in \{1, \cdots, N\}$ one
 can find $u^\prime_i \in C^0(\overline{V_i} , \widehat{E})$
 such that  $\forall \theta \in U_i^\prime$ one has:
 $$
 (d)\;\; |u_i(\theta, \cdot)|_{L^2}=1,\; |\chi_1(D)(u_i(\theta, \cdot))|_{L^2}
 > \frac{2}{3},\; |u_i(\theta, \cdot)|_{H^{-1}} < \epsilon.$$
 Now consider for each $i \in \{1, \cdots, N\}$
 $\phi_i \in C^\infty_c(U_i^\prime)$ such that $\phi_i \equiv 1$ on
 $U_i.$ Set
 $$
 u \,=\,\frac{\sum_{i=1}^N \phi_i u_i}{\sqrt{\sum_{i=1}^N \phi^2_i}}.
 $$ Then properties (a), (b), (c), (d) and the proof of
 Lemma 4 of \cite{LPCUT} show that if $\epsilon >0$ is small enough
 then $\langle \chi_1(\D)(u) ;\chi_1(\D)(u) \rangle$ is invertible in $C(T)\rtimes_r \Gamma$
 which proves the result.
\end{proof}

\subsection{Difference class associated to two spectral sections}

\begin{proposition}\label{max-min} (Wu)
Let $\P_1$, $\P_2$ be two spectral sections for $\D$. Then there exists
  spectral section $\Q$, $\R$ for  $\D$ such that for $j\in\{1,2\}$:
$\P_j \R = \P_j = \R \P_j$  and
$\P_j \Q  = \Q =  \Q \P_j$.
\end{proposition}

\begin{proof} We follow closely the unpublished proof of Wu (\cite{Wu}).
\begin{lemma} \label{Fan} We set $\E= L^2_{\Cr}(\wm,\we).$ Assume the existence of a
spectral section for $\D$.
{\item 1)}  For any spectral cut
$\chi_1$, there is a smooth spectral cut $\chi_2$ with
$\chi_1(t) \chi_2(t)=\chi_1(t)$, and a spectral section $\R$
satisfying:
$$
\chi_1(\D) (\E) \subset \R(\E) \subset \chi_2(\D)(\E).
$$
{\item 2)} Similarly, for any smooth spectral cut $\chi_2$, there a
spectral cut $\chi_1$ with $ \chi_1(t) \chi_2(t)=\chi_1(t)$ and a
spectral section $\Q$ satisfying
$$\chi_1(\D) (\E) \subset \Q(\E) \subset \chi_2(\D)(\E).
$$
\end{lemma}
\begin{proof} 1) Let $\P$ be a spectral section for $\D$ satisfying
$$ g_1(\D) (\E) \subset \P(\E) \subset g_2(\D)(\E) $$ where $g_1,
g_2$ are smooth spectral cuts. Let $\chi_1$ be a smooth spectral
cut. Choose a smooth spectral cut $\chi$ satisfying: $$ \chi(t)
\chi_1(t) = \chi_1(t)\; {\rm and} \; \chi(t) g_2(t) = g_2(t). $$
We have $\chi(\D) \P=\P \chi(D)=\P.$ 
Thus: \begin{align*} (\id -\chi(\D))(\id
-\P)&=(\id-\P)(\id-\chi(D))=\id-\chi(\D) \\&=(\id-\P)(\id-\chi(\D))(\id-\P).
\end{align*}
Working in the $C^*-$algebras $\B_{\Cr}( ( \id-\P)(\E))$ and
$$\K_{\Cr}( ( \id-\P)(\E))=(\id -\P) \K_{\Cr}(\E)(\id -\P),$$ let
$\{P_n\}\subset (\id -\P) \K_{\Cr}(\E)(\id -\P)$ be an approximate
unit, $\P_n \leq \P_{n+1}$ and $\P_n$ are projections in
$\K_{\Cr}( ( \id-\P)(\E))$. Then we have $$( (\id-P)-(\id
-\chi(\D))\P_n \rightarrow (\id-\P)-(\id -\chi(D))$$ in norm in
$(\id -\P) \K_{\Cr}(\E)(\id -\P).$ Let $N_0$ be such that $$ ||
((\id -\P) -P_{N_0})-(\id-\chi(\D))((\id -\P) -P_{N_0})|| <
\frac{1}{2}. $$ As $(\id-\P)-\P_{N_0}$ is a projection in
$\B_{\Cr}((\id -\P)(\E),$ Lemma 5 of \cite{LPCUT} (due to Wu)
implies that there is a projection $\R_0$ in  $\B_{\Cr}((\id
-\P)(\E))$, $\R_0-(\id-\P) \in \K_{\Cr}((\id -\P)(\E))$, such that
\begin{align*}(\id-\R_0)((\id-\P)(\E))&=(\id-\chi(\D))((\id -\P)-\P_{N_0})(\id
-\P)(\E) \\&\subset (\id -\chi(\D))(\E). 
\end{align*}
Let $\R_1\in
\B_{\Cr}(\E)$ be the projection which is $\R_0$ on $(\id-\P)(\E)$
and $\id$ on $\P(\E).$ Then $\R_1-\P\in \K_{\Cr}(\E)$ and $$ \R_1
\chi_1(\D) = \chi_1(\D) \R_1=\chi_1(\D) \;{\rm and}\; \R_1
g_2(\D)=g_2(\D) \R_1=g_2(\D). $$ To modify $\R_1$ to the desired
spectral section, let $\psi_N(t):= \psi(t/N)$ where $\psi$ is a
smooth spectral cut such that $\psi(t)\equiv 1$ on $[-1,
+\infty[.$ Then since $\D$ is a regular operator, one checks
easily that as $N\rightarrow +\infty$ $$ \psi_N(\D
)(\chi_1(\D)-\R_1) \rightarrow \chi_1(\D)-\R_1,\;\, \psi_N(\D
)(g_2(\D)-\R_1) \rightarrow g_2(\D)-\R_1 $$ in norm in
$\K_{\Cr}(\E).$ Let $N_0$ be chosen such that $\psi_{N_0}
\chi_1=\chi_1$ and $\psi_{N_0} g_2=g_2$ and $$ ||\R_1-\chi_1(\D) +
\psi_{N_0}(\D)( \chi_1(\D)-\R_1)||=||R_1-\psi_{N_0}(\D) \R_1||<
\frac{1}{2}. $$ Applying Lemma 5 of \cite{LPCUT}, we get a
projection $\R=\R^2=\R^*\in \B_{\Cr}(\E)$, $\R-\R_1\in
\K_{\Cr}(\E),$ such that $$ \R(\E)=\psi_{N_0}\R_1(\E)\subset
\psi_{N_0}(\D)(\E). $$ Let $\chi_2(t):=\psi_{N_0}(t).$ Then
$\chi_2$ is a smooth spectral cut and we have $\R(\chi)\subset
\chi_2(\D)(\E).$ On the other hand, since $$
\chi_1(\D)=\psi_{N_0}(\D)\chi_1(\D)=\psi_{N_0} \R_1 \chi_1(\D), $$
we get $\chi_1(\D)(\E) \subset \psi_{N_0}(\D) \R_1(\E)=\R(\E).$
Similarly, we also have $g_2(\D)(\E) \subset \R(\E).$ Therefore,
$\R$ is a spectral section with the desired property: $$\chi_1(\D)
(\E) \subset \R(\E) \subset \chi_2(\D)(\E). $$ 2) is proved
similarly
\end{proof}

We go back to the proof of Proposition \ref{max-min}. There are
smooth spectral cuts $g_1,\,g_2$ such that: $$ g_1(\D)(\E) \subset
\P_j(\E) \subset g_2(\D)(\E),\; j=1,2. $$ Applying Lemma
\ref{Fan}, we find smooth spectral cuts $\chi_1,\chi_2$ with
$\chi_1 \cdot g_1=\chi_1$, $\chi_2 \cdot g_2=g_2$ and spectral
sections $\R,\Q$ with $$ \chi_1(\D)(\E)\subset \Q(\E)\subset
\P_j(\E) \subset g_2(\D)(\E)\R(\E)\subset \chi_2(\D)(\E). $$ Now
one checks easily that $\Q,\R$ satisfy the desired property:
$$\P_j \R = \P_j = \R \P_j,\; \P_j \Q  = \Q =  \Q \P_j. $$

\end{proof}

Recall  the following stability result:
$$
K_0( \K_{\Cr} (\E) ) \simeq K_0(\Cr).
$$

\begin{definition} Let $\P_1$ and $\P_2$ spectral sections for $\D.$
Then there exists a difference class
$[\P_1 -\P_2] \in K_0(\Cr)$ defined as follows. Choose a
spectral section $\P^\prime$ of $\D$ such that
$\P_i \P^\prime = \P^\prime = \P^\prime \P_i$ for
$i \in \{1,2\}$. Then $\P_1 -\P^\prime$ and $\P_2 -\P^\prime$ induce
projections in $\K_{\Cr}$ and
$[\P_1 -\P_2]= [\P_1 -\P^\prime] -[\P_2 -\P^\prime]$ is well defined as
an element of $K_0(\Cr).$
\end{definition}

\subsection{The noncommutative spectral flow on foliated bundles.}

Now we consider a continuous family $(g_u)_{u\in [0,1]}$ of
$\Gamma-$equivariant vertical metrics of the fibration
$\widehat{M} \rightarrow T.$ We consider also a continuous family
$(h_u)_{u\in [0,1]}$ of $\Gamma-$equivariant hermitian metrics on
the Clifford module $\we \rightarrow \wm$ and a continuous family
of $\Gamma-$equivariant hermitian connections $(\nabla_u)_{u\in [0,1]}$ on $\we$. For each
$u\in[0,1]$ one then gets a $\Cr-$Hilbert module
 $L^2_{\Cr}(\wm , \we,u)$ which depends on $u$ via $g_u$ and
 $h_u$. These Hilbert modules are all isomorphic, being
each one isomorphic to the standard Hilbert $\Cr$-module
$\mathcal{H}_{\Cr}$.
 Let
$(\D_u)_{u \in [0,1]}$ be the associated  family of
$\Cr-$linear Dirac type operators where each $\D_u$ acts on
$L^2_{\Cr}(\wm , \we,u).$ This family is continuous in the following sense.

\smallskip
Let $\{U_j, \; 1\leq j \leq l\}$ be a finite set of open subsets of
$\wm$ satisfying the three following properties:

(i) Each $U_j$ is diffeomorphic to the open ball $B(0,1)$ of $\RR^n.$

(ii) $\cup_{ 1\leq j \leq l} U_j \cdot \Gamma = \wm$

(iii) For each $j \in \{1, \ldots, l\}$, the restriction to $U_j$ of the
following bundles is trivial: $T Z \rightarrow \wm$, $\we \rightarrow \wm$,
$T \wm \rightarrow \wm$.

Then the restriction of $\D_u$ to each
$U_j$ induces a differential operator of order one
acting on $C^\infty(B(0,1) ; \CC^N)$ (for a suitable $N$): $D_{j,u}
=\sum_{k=1}^n a_k(z,u) \partial_{z_k} + b(z,u)$ ($z \in B(0,1)$).

Then we observe  that the coefficients of  $D_{j,u},$
$u\rightarrow a_k(z,u)$, $u\rightarrow b(z,u)$ belong
to $C^0( [0,1] ; C^\infty(B(0,1) ; \CC^N) ).$ We then say that
the family $(\D_u)_{u \in [0,1]}$ is continuous in this sense (the definition
does not depend on the choice of the $U_j$ and of the trivializations).

\smallskip
 We assume that the
index class defined by one (and thus any) of the $\D_u$ is
trivial. Recall that the family $(\D_u)_{u \in
[0,1]}$ defines in a standard way a $(C^0 [0,1]\otimes \Cr)-$linear
operator acting on the $(C^0 [0,1]\otimes \Cr)-$Hilbert module 
defined by the bundle $\mathcal{E}$ over $[0,1]$ with fiber 
$\mathcal{E}_u=L^2_{\Cr}(\wm , \we,u),\; u \in [0,1].$ 
This operator is of Dirac-type and therefore  $(C^0 [0,1]\otimes \Cr)-$Fredholm
by the Mishchenko-Fomenko calculus.
Recall also that 
we have a natural isomorphism
$$
K_1 (C^0 [0,1]\otimes \Cr) \simeq K_1 ( \Cr)
$$ 
which is implemented by the evaluation map
$f(\cdot)\otimes\lambda\to f(0)\lambda$.
This implies that the total index class of $(\D_u)_{u \in
[0,1]}$, in $K_1 (C^0 [0,1]\otimes \Cr)$, is also zero. By the
existence theorem this implies that the family $(\D_u)_{u\in
[0,1]}$ admits a (total) spectral section $(\P_u)_{u\in
[0,1]}.$ The following definition is inspired by the work of Dai-Zhang
\cite{Dai-Zhang}:

\begin{definition}
If $\Q_0$ (resp. $\Q_1$) is a spectral section associated with
$\D_0$ (resp. $\D_1$) then the noncommutative spectral flow 
from $(\D_0,\Q_0)$ to
$(\D_1,\Q_1)$ through $(\D_u)_{u\in [0,1]}$, is the
$K_0(\Cr)-$class: $$ {\rm sf}( (\D_u)_{u\in[0,1]}; \Q_0, \Q_1)=
[\Q_1 -\P_1] - [\Q_0-\P_0] \in K_0(\Cr). $$
\end{definition}

Proceeding as in the work of Dai-Zhang one can prove that the definition
does not depend on the choice of total spectral section $(\P_u)_{u\in
[0,1]}.$

\subsection{Trivializing perturbations}\label{subsect:perturbations}

Let $\Ind(\D)$ be equal to zero.
The operator $\D$ will not be, in general, invertible.
 Let $\P$ be a spectral section
for $\D$; then $\P$ fixes
a specific trivialization of $\Ind(\D)$; this is achieved by
defining a perturbation $\A^0_{\P}$ of $\D$ such that
$\D+\A^0_{\P}$ be invertible. This subsection is devoted
to make this statement precise. First we introduce the
relevant space of pseudodifferential operators.

\begin{definition}\label{def:pseudos}
{\item 1)} Denote by $\Psi^{-\infty}_{\Cr}( \wm , \we)$ the set of
operators $R \in \B_{\Cr}$ such that for any $N\in \NN$, $R$
extends as a continuous operator $$H^{-N}_{\Cr}(\wm , \we) 
\longrightarrow H^{N}_{\Cr}(\wm , \we)\,.$$ {\item 2)} Let $k \in \ZZ.$ Denote
by $\Psi^k_{\Cr}(\wm , \we)$ the set of bounded operators $A:
H^k_{\Cr}(\wm , \we) \rightarrow L^2_{\Cr}(\wm , \we)$ satifying
the following property. There exists $\epsilon >0$ and one can
write $A=B + R$ where $R \in \Psi^{-\infty}_{\Cr}( \wm , \we)$ and
$B=(B_\theta)_{\theta \in T}$ is a smooth $\Gamma-$equivariant
family of fiberwise pseudo-differential operators of order k such
that the Schwartz kernels of each $B_\theta$ vanish outside an
$\epsilon-$neighborhood of the diagonal.
\end{definition}

\s
\n
{\bf Remarks.}

\n
1) The operator $B$ appearing in 
definition \ref{def:pseudos} is an element in  $\Psi^k_{\Cr,c}(\wm , \we)$,
the space of $\Gamma$-equivariant families of order $k$
pseudodifferential operators with Schwartz kernel of {\it compact
$\Gamma$-support}, i.e. the support of the Schwartz kernel,
viewed as an element in $\wm\times_\pi \wm$
($\pi$ denotes the projection $\wm\to T$),
defines a compact
set  in $\wm\times_\pi \wm/\Gamma$.

\n
2) One can prove, see for example Morioshi and Natsume \cite[Section 3]{Mor},
that a $\Gamma$-equivariant family of Dirac operators
$(D(\theta))_{\theta\in T}$ admits a parametrix
$(Q(\theta))_{\theta\in T}$ in $ \Psi^{-1}_{\Cr,c}(\wm , \we)$
with rests in $\Psi^{-\infty}_{\Cr,c}(\wm , \we).$  Moreover,
an element $S=(S(\theta))_{\theta\in T}\in \Psi^k_{\Cr,c}(\wm , \we)$
defines a bounded
$\Cr-$linear operator from $H^m_{\Cr}(\wm,\we)$ to  $H^{m-k}_{\Cr}(\wm,\we)$.
It is precisely  this result that is used in order to define
the index class associated to $(D(\theta))_{\theta\in T}$.

\s
Proceeding as in the proof of Propositions 2.5 and 2.10 of
Leichtnam-Piazza \cite{LPGAFA} (and thus, ultimately,
as in Lemma 8 of \cite{MP I}) one can prove the following

\begin{proposition} \label{Triv} Assume that $\Ind \D=0$ in $K_1(\Cr)$ and
let $\P$ be a spectral section for $\D$. Then $\P \in
\Psi^0_{\Cr}(\wm, \we)$ and there exists a self-adjoint operator
$\A^0_{\P} \in \Psi^{-\infty}_{\Cr}(\wm,\we)$ with the following
three properties. {\item (1)} We can find a real $R>0$ such that
$\varphi(\D)\circ \A^0_{\P}\equiv 0$ for any function $\varphi \in
C^\infty(\RR ,\CC)$ vanishing on $[-R , R].$ {\item (2)} $ \D +
\A^0_{\P}$ is invertible from  $H^1_{\Cr}(\wm , \we)$ onto
$L^2_{\Cr}(\wm , \we)$.
 {\item (3)} $\P$ is equal to
the projection onto the positive part of $ \D + \A^0_{\P}$, $$ \P=
\frac{1}{2}\left( \frac{\D + \A^0_{\P}}{|\D + \A^0_{\P}|} +\id
\right). $$
\end{proposition}

We shall not reproduce the proof of this proposition here.
We simply add  that the construction of $\A^0_{\P}$
does depend on choices; however, if $\mathcal{B}^0_{\P}$
is a different trivializing perturbation associated to
$\P$, then \begin{equation}\label{independence}
\text{the operator}\;\;\D+r \A^0_{\P}+ (1-r)\mathcal{B}^0_{\P}\;\;
\text{is invertible}\;\;\forall r\in [0,1].
\end{equation}

\section{Index classes on foliated bundles with
boundary}\label{sect:b-aps-classes}

\subsection{Preliminaries: numerical indeces on manifolds with boundary.}
\label{subsect:preliminaries}$\;$

\smallskip

{\bf The invertible case}. In order to orient the reader we review in some detail
the various indeces that can be attached to a Dirac operator
on an even dimensional manifold with boundary.

We thus consider
 a smooth connected compact manifold with boundary $M.$
We fix a boundary defining function
$x\in C^\infty(M)$. Let
$g_M$ be a riemannian metric on $M$; we assume this metric
to be of product type near the boundary. We consider
on $M$ a
 unitary Clifford module $E$ endowed with a unitary
connection $\nabla^E$ which is Clifford with respect to
the Levi-Civita connection associated to the metric $g_M$.
We obtain in this way a Dirac-type operator $D$.
Suppose now that $M$ is {\it even} dimensional
so that $E$ is $\ZZ_2$-graded: $E=E^+\oplus E^-$.
The Clifford bundle associated to $T^*(\pa M)$
and to the boundary metric
acts in a natural way on $E_{|\pa M}$: $$\forall e \in E_{|\pa M},\;\;
\forall \eta \in T^*\pa M,\;\;
\cl_{\pa}(\eta)(e_{|\pa M}):= \cl(
dx)\cl(\eta)(e_{|\pa M})$$
 We define $E_0$ to be $E_{|\pa M}^+$. It is a
unitary Clifford bundle with respect to $\cl_{\pa}(\cdot)$.
It is endowed with the  induced Clifford connection.
We denote  by $D_0$ the associated Dirac operator and we
call it {\it the boundary operator of $D.$}
Finally, we identify $E_{|\pa M}^-$ with
$E_0$ through Clifford multiplication by
$ \cl(i dx)$, denoted in the sequel by $\sigma$.
With these identifications the operator $D^+$ can be
written near the boundary as $\sigma(\partial_x + D_0)$.

\m
As already remarked the APS-boundary value problem is
obtained by considering
the operator $D^+$ with domain $$\{u\in C^\infty (M,E^+)\;|\;
u_{|\pa M}\in {\rm Ker} \Pi_{\geq}\}$$ with $\Pi_{\geq}=\chi_{[0,\infty)} (D_0)$.
Let $\ind(D^+,\Pi_{\geq})$ be the APS-index and assume for the time being
that $D_0$ is {\it invertible}.
Then we can describe this index in a different way:
we can attach an infinite cylinder $(-\infty,0]\times \partial M$ to $M$
along its boundary $\partial M$, thus obtaining a manifold with
cylindrical ends $M_{{\rm cyl}} $ and with product metric $dx^2 +
g_{\partial M}$ along the cylinder. The operators $D$ extends
in a natural way to an operator $D_{{\rm cyl}}$
on the manifold $M_{{\rm cyl}}$, acting on the sections
of the bundle $E_{{\rm cyl}}$ obtained by extending in
an obvious way $E$. It turns out
that this operator is Fredholm, as a bounded linear map from
$H^1 (M_{{\rm cyl}},E_{{\rm cyl}})$ to $L^2 (M_{{\rm cyl}},E_{{\rm cyl}}) $,
and that its index is equal to the APS-index:
\begin{equation}
\ind (D^+,\Pi_{\geq})=\ind (D_{{\rm cyl}}^+)\,.
\end{equation}
 This equality
is explained in the original paper of Atiyah-Patodi-Singer
where it is proved, more precisely, that the kernel and cokernel of
the two operators are naturally isomorphic.
Connected with the cylindrical picture is Melrose' $b$-picture
\cite{Melrose}:
the change of
coordinates $x=\log y $ compactifies $ M_{{\rm cyl}}$ to a
compact manifold
with boundary ${}^b M$ but with a degenerate metric ${}^b g$
which can be written as $dy^2/y^2 + g_{\partial
M}$ near the boundary. The operator on $M_{{\rm cyl}} $ then defines
in a natural way a differential operator ${}^b D$
on the compactified
manifold ${}^b M$; up to a bundle isomorphism
the operators ${}^b D^\pm$ can be written near the boundary
as $\pm y\pa_y +D_0$; this means that  ${}^b D$ is
generated by the vector fields on ${}^b M$ which are tangent to the boundary;
${}^b D$ is therefore, by definition, a $b$-differential operator.
Melrose has developed on ${}^b M$
a pseudodifferential calculus,
 which extends the algebra of $b$-differential
operators; this is known as the $b$-calculus and it can be used,
among other things,
in order to show that ${}^b D$
is Fredholm on naturally defined $b$-Sobolev spaces
$H^m_b ({}^b M, {}^b E)$,
with index
equal to the APS-index.
It should certainly be remarked that
the cylindrical picture and the $b$-picture
are two different descriptions of the same mathematical
object.
\\
{\it Summarizing}:
if the boundary operator $D_0$ is invertible
\begin{equation}\label{equality-invertible}
\ind(D^+,\Pi_{\geq})=\ind({}^b D^+)=\ind (D_{{\rm cyl},+})\,.
\end{equation}

\m

{\bf The general case.}
In the $b$-picture ($\equiv$ cylindrical picture)
it is fundamental to assume that $D_0$ is invertible; if the kernel
of $D_0$ is non-trivial, then ${}^b D$ will {\it not} define a Fredholm
operator.
Still, the APS-index is indeed equal to a $b$-index but for
 a {\it perturbed} $b$-operator;
we shall now describe this fundamental point
in full generality, thus considering the APS-index
$\ind(D^+,P)$ associated to an arbitrary
 spectral section $P$ for $D_0$. \\As already remarked in Subsection
\ref{subsect:perturbations}
one  can prove that there is
a smoothing operator $A^0_P$ on the boundary $\pa M$ such that
$D_0+A^0_P$ is invertible. The perturbation $A^0_P$ can be extended
from the boundary to the interior, thus defining a smoothing
$b$-operator $A_P^+$. The construction of this operator will
be recalled below, in the general case of foliated bundles.
 The operator ${}^b D^+_P :=
{}^b D^+ + A_P^+$ is now Fredholm
on $b$-Sobolev spaces $H^m_b$, with  index independent of $m$
and equal to the APS-index
$\ind(D^+,P)$. See \cite{MP I} for proofs and details.\\
As an example, consider $P=\Pi_{\geq}$ but assume that
 ${\rm Ker} D_0\not= 0$; then $A^0_P$ is nothing but the $L^2$-orthogonal
projection onto ${\rm Ker} D_0$.

 Notice that we could extend
$A^0_P$ to an operator on the cylindrical manifold $M_{{\rm cyl}}$
by simply employing  a cut-off function $\phi$
 equal to 1 on the attached half-cylinder and
equal to zero on the complement of the collar neighbourhood
of the boundary of $M$.
The resulting operator $D_{{\rm cyl},+}+ \sigma A^0_P \phi$
can be either viewed as an operator on the manifold $M_{{\rm cyl}}$
or as an operator on ${}^b M$; as such
it does not define
a $b$-pseudodifferential operator; however, it is still
possible to prove that it is Fredholm as a map
$H^1 (M_{{\rm cyl}},E_{{\rm cyl}}) \to L^2 (M_{{\rm cyl}},E_{{\rm cyl}})$ or,
equivalently, as a map  $H^1_b ({}^b M, {}^b E)\to L^2_b ({}^b M, {}^b E)$
and with index equal to
$\ind(D^+,P)$.
\\
{\it Summarizing:}\\
if $P$ is a spectral section for the boundary operator $D_0$ then
\begin{equation}\label{equality}
\ind(D^+,P)=\ind({}^b D^+ + A_P^+)=\ind (D_{{\rm cyl},+}+ \sigma A^0_P \phi)\,.
\end{equation}

We remark that in order to establish an {\it
index
formula} for one of these 3 indeces, it is  very useful to consider
the $b$-perturbation $A^+_P$; such a formula
is obtained
in \cite{MP I}.

\subsection{Foliated bundles with boundary.}\label{subsect:foliated-with-boundary}

Let $\Gamma$ be a finitely generated discrete group. Let $T$ be a
smooth closed compact connected manifold on which $\Gamma$ acts on
the right. Let $\widehat{M}$ be a manifold {\it with boundary} on
which $\Gamma$ acts freely, properly and cocompactly on the right:
the quotient space $M=\widehat{M} /\Gamma$ is thus a smooth
compact manifold with boundary. We assume that $\widehat{M}$
fibers over $T$ and that the resulting fibration $$ \pi:
\,\widehat{M} \rightarrow T $$ is a {\it $\Gamma$-equivariant
fibration} with fibers $\pi^{-1}(\theta),\theta \in T,$ that are
transverse to $\partial \widehat{M}$ and of dimension $2k$
(on a manifold with boundary we are always assuming
that the fibers are even-dimensional). Notice
that each fiber is a smooth manifold with boundary; we shall also
denote the typical fiber of $\pi: \,\widehat{M} \rightarrow T $ by
$Z$.
We choose a $\Gamma$-invariant product-like metric
on the vertical tangent bundle $TZ $.
Finally, we  assume the existence of a $\Gamma-$equivariant spin
structure on $TZ $ that is fixed once and for all. We denote
by $S^Z \rightarrow \widehat{M} $ the associated spinor bundle.


\smallskip
The compact manifold with boundary $M$ inherits a foliation
${\cal F}$, with leaves equal to the image of the fibres
of $\pi: \,\widehat{M} \rightarrow T
$ under the quotient map $\widehat{M}\rightarrow M=\widehat{M}/\Gamma$.
Notice that the foliation ${\cal F}$ is transverse to the
boundary of $M$.

\smallskip \noindent {\bf Example.} Let $X$ be a compact manifold
with boundary and let $\Gamma\rightarrow \widetilde{X}\rightarrow
X$ be a Galois cover of $X$. Let $T$ be a smooth compact manifold
on which $\Gamma$ acts by diffeomorphisms. We consider
$\widehat{M}=\widetilde{X}\times T$, $ \pi=\text{projection onto
the second factor}$, $M=\widetilde{X}\times_{\Gamma}
T:=(\widetilde{X}\times T)/\Gamma$ where we let $\Gamma$ act on
$\widetilde{X}\times T$ diagonally. As a particular example of
this construction consider a smooth closed riemann surface
$\Sigma$ of genus $g>1$ and let $\Gamma=\pi_1 (\Sigma)$, a
discrete subgroup of $PSL(2,\RR)$. Let $\{p_1,\dots,p_k\}$ points
in $\Sigma$ and let $D_j\subset \Sigma$ be  a small open disc
around $p_j$. Let $D=\cup_{j=1}^{k} D_j$. Then we can consider
$X:=\Sigma\setminus D$, $\Gamma\rightarrow\widetilde{X}
\rightarrow X$ the Galois cover induced by the universal cover
$\HH^2\rightarrow \Sigma$, $T=S^1$, with $\Gamma$ acting on $S^1$
by fractional linear transformations.

\smallskip
Next we consider a $\Gamma-$invariant boundary defining function
$x$ of $\partial \wm$ and a $\Gamma$-equivariant complex hermitian
vector bundle $\widehat{V} \rightarrow \widehat{M} $ endowed with
a $\Gamma-$invariant hermitian connection $\widehat{\nabla}.$
We then set $\widehat{E}= S^Z \otimes
\widehat{V}= \widehat{E}^+ \oplus \widehat{E}^-$ which defines a
smooth $\Gamma-$equivariant family of  $\ZZ_2-$graded hermitian
Clifford modules on the fibers $\pi^{-1}(\theta), \, \theta \in
T$. We then get a smooth family of $\Gamma-$equivariant
$\ZZ_2-$graded Dirac type operators $$  D(\theta)= \begin{pmatrix}
0 & {D^-(\theta)} \cr {D^+(\theta)}  & 0 \cr
\end{pmatrix},\; \theta\in T
$$
 acting fiberwise on
$C^\infty_c( \widehat{M},\, \widehat{E})$.
 Moreover in a collar neighborhood ($ \simeq [0,1] \times \partial \pi^{-1}(\theta) =
\{(x, y)\}$) of  $\partial \pi^{-1}(\theta) $ we may write:
$$
 {D^+(\theta)}= \sigma (\partial_x + D_0(\theta)\,)
$$ where $ D_0(\theta)$ is the induced boundary Dirac type operator acting
on
$$
 C^\infty(\partial \pi^{-1}(\theta) ,\, \widehat{E}^+_{|_{\partial \pi^{-1}(\theta)}}\,).
$$ Observe that our family can also be thought as a longitudinal
operator on $(M,{\cal F})$ acting on the sections of
$E:=\widehat{E}/\Gamma$.

\smallskip The family $(D(\theta))_{\theta \in T}$ defines a
$\Cr-$linear $\ZZ_2-$graded Dirac-type operator $\D$ acting
on $C^\infty_c( \wm , \we).$ Similarly,
the family $(D_0(\theta))_{\theta \in T}$ defines a
$\Cr-$linear Dirac type operator $\D_0$ acting
on the space $C^\infty_c(\partial \wm , \we^+_{|\partial \wm}).$

\m
Corresponding to the above metrically-incomplete picture, there is
a $b$-picture, obtained by attaching an infinite cylinder $(-\infty,0]\times \pa \wm$
to $\widehat{M}$
and compactifying it as we did in the previous subsection. We shall keep
the same notation for the resulting manifold; we shall denote by
${}^b \D$ the $C(T)\rtimes_r \Gamma$-linear operator
defined by the $\Gamma$-equivariant family of $b$-Dirac operators
$${}^b D (\theta)=\sigma(y\partial_y +D_0 (\theta))$$
associated to the  $b$-data.

In our recent paper \cite{LPETALE} we develop a $b$-calculus $\Psi^*_{b,\Cr}$
 on foliated bundles
with boundary; we then employ such a calculus in order to establish
the following

\begin{theorem}
Assume that there exists a real $\epsilon>0$ such that
for any $\theta\in T$, the $L^2$-spectrum
of $D_0 (\theta)$ acting on
$L^2(\partial\pi^{-1}(\theta),\widehat{E}^+_{|_{\partial\pi^{-1}(\theta)}})$
does not meet $]-\epsilon,\epsilon[$.
Then $({}^b D^+(\theta))$ defines
$\forall m \in \NN^*$ a $C(T) \rtimes_r \Gamma$-linear bounded
operator $${}^b  \D: H^m_{b, C^0(T) \rtimes_r \Gamma}(\wm,\widehat{E}^+)
\longrightarrow  H^{m-1}_{b, C^0(T) \rtimes_r
\Gamma}(\wm,\widehat{E}^-)$$ which is invertible modulo $C^0(T)
\rtimes_r \Gamma $-compacts.
There is a well defined  $b$-index class $\Ind_b ({}^b \D^+)$ in
$K_0(C^0(T) \rtimes_r \Gamma)$, independent of $m$.
\end{theorem}

 Let us consider the operator $\D_0$
defined by the boundary family $(D_0 (\theta))_{\theta\in T}$;
it defines a (regular) unbounded operator
on the $C(T)\rtimes_r \Gamma$-Hilbert module
$L^2_{C(T)\rtimes_r \Gamma}$.
One can prove (see \cite[Proposition 1]{LPETALE}) that if
the hypothesis of the above theorem  holds then $\D_0$
is $L^2_{C(T)\rtimes_r \Gamma}$-invertible with
domain $H^1_{C(T)\rtimes_r \Gamma}$ and that its
inverse is  induced by the $\Gamma$-equivariant family of operators
$\{D_0 (\theta))^{-1}\}_{\theta\in T}$.
 We can thus consider
the $C(T)\rtimes_r \Gamma$-linear bounded operator
\begin{equation}\label{aps-invertible}
\Pi_{\geq} (\D_0):= \frac{1}{2} \left( \frac{\D_0}{|\D_0|}+\Id
\right)\,.
\end{equation}
This is in fact a self-adjoint projection and can be used in order
to define an APS-index class  $\Ind^{APS} (\D^+,\Pi_{\geq} (\D_0))\in K_0 (\Cr)$;
we shall see the details below. We shall also see that
$$
\Ind_b ({}^b \D^+)=\Ind^{APS} (\D^+,\Pi_{\geq} (\D_0))\in K_0 (\Cr)\,.$$

Thus in the {\it invertible case} we can extend to the
present  noncommutative context
 the basic results recalled for the numerical indeces
 in Subsection \ref{subsect:preliminaries}.
In the non-invertible case the operator
${}^b  \D: H^m_{b, C^0(T) \rtimes_r \Gamma}(\wm,\widehat{E}^+)
\longrightarrow  H^{m-1}_{b, C^0(T) \rtimes_r
\Gamma}(\wm,\widehat{E}^-)$ will {\it not} be $\Cr$-Fredholm; similarly,
as already remarked, in the non-invertible case
the operator (\ref{aps-invertible}) does not make sense as
a bounded $\Cr$-linear operator and we cannot define a APS-index class.
The way out in the general case is therefore to consider spectral sections
$\P$
for the boundary operator $\D_0$; the existence of these spectral sections
follow from our basic result, Theorem \ref{theo:existence},
and the cobordism invariance
of the index class associated to $\D_0$ in $K_1 (\Cr)$, a result that will be established
in the next subsection.

\subsection{Cobordism invariance}\label{subsect:cobordism}

\begin{theorem}\label{theo:cobordism}
Let $\D_0$ be the $\Cr$-linear operator defined
by the boundary family $(D_0 (\theta))_{\theta\in T}.$
One has $\Ind \,\D_0=0$ in $K_1(\Cr).$
\end{theorem}

\begin{proof} The proof employs equivariant $KK$-theory.
It is easy to see that
$\D_0$ defines a class $[\D_0]$ in the $\Gamma-$equivariant
Kasparov
group $KK_\Gamma^1(C_0(\pa \wm) , C(T)).$ Recall  that, since
$C_0(\pa \wm) \rtimes \Gamma$ is Morita equivalent
to $C(\pa M),$ one  has a
natural map $\Theta: KK_\Gamma^1(C_0(\pa \wm) , C(T)) \rightarrow
KK^1(C(\pa M) , \Cr).$
If
$\pi^{\pa M}: \pa M \rightarrow $pt denotes the mapping of $\pa M$
to a point,  then, under the natural isomorphism $KK^1(\CC,
\Cr)\simeq K_1(\Cr),$ we have $\Ind \D_0=\pi_*^{\pa M}\circ \Theta([\D_0]).$

Let
$C_{0,\pa \wm}(\wm)\subset C_0(\wm)$ denote the ideal of continuous
functions on $\wm$ vanishing on the boundary, let $i$ be the natural
inclusion of $\pa \wm$ into $\wm$ and consider the long exact
sequence, in $KK_\Gamma( \cdot , C(T))$, associated to the semisplit short
exact sequence:
\begin{equation} \label{ext}
0 \rightarrow C_{0,\pa \wm}(\wm)\overset{j}\rightarrow  C_0(\wm)
 \overset{q} \rightarrow C_0(\pa
\wm)\rightarrow 0
\end{equation} (see Blackadar \cite{Blackadar} page 197 and Chapter 20). We have in
particular the exactness of
$$ KK_\Gamma^0(C_{0,\pa \wm}(\wm),
C(T))\overset{\delta_\Gamma}\mapsto  KK_\Gamma^1(C_0(\pa \wm), C(T)) $$ $$
\overset{\iota_*}\mapsto KK_\Gamma^1(C_0( \wm), C(T)) $$ and thus $i_*\circ
\delta_\Gamma=0.$
\begin{lemma} We have $[\D_0]=\delta_\Gamma[\D]$ where
$[D]\in KK_\Gamma^0(C_{0,\pa \wm}(\wm), C(T))$ is the class defined by $\D$.
\end{lemma}
\begin{proof} We are using  both a $\Gamma-$equivariant and bivariant
generalization of the proof of Theorem 5.1 of Higson \cite{Higson}. We can
replace $\wm$ by a $\Gamma-$equivariant collar neighborhood $\widehat{W}$ 
($\simeq
[0,1) \times \pa \wm$)
of $\pa \wm$  such that the restriction of $\pi$ to $\widehat{W}$
induces a $\Gamma-$equivariant fibration over $T.$ Consider the differential
operator $d$
$$ d=\begin{pmatrix}
 0 & - i \frac{d}{d x} \cr
-i \frac{d}{d x}& \cr
\end{pmatrix}
$$ acting on $[0,1]$. It defines a class in $KK_\Gamma ^1( C_0(0,1), \CC)$.
Recall that
the Kasparov product $[d]\otimes \cdot$ induces an isomorphim:
$$
[d]\otimes: KK_\Gamma^1(C_0(\pa \wm), C(T))\mapsto
KK_\Gamma^0(C_{0,\pa \wm}(\widehat{W}), C(T))
$$
As in \cite{Higson},
the connecting map $\delta_\Gamma:$
$$
KK_\Gamma^0(C_{0,\pa \wm}(\widehat{W}),
C(T))\overset{\delta_\Gamma}\mapsto  KK_\Gamma^1(C_0(\pa \wm), C(T))
$$ is given by the inverse of $[d]\otimes \cdot.$
Denote by $\D_{\widehat{W}}$ the restriction of $\D$ to
$\widehat{W}$, then  one checks (as in Theorem 4.7 of \cite{Higson})
that $\D_{\widehat{W}}= [d]\otimes [\D_0]. $ One then
gets $\delta_\Gamma [ \D_{\widehat{W}}]= [\D_0]$ which proves the result.
\end{proof}

 We also have a natural map, still denoted $\Theta:$
$$
\Theta: KK_\Gamma^1(C_0( \wm), C(T))\rightarrow KK^1(C( M),
C(T)\rtimes \Gamma).
$$
Denote by $\pi^M$ the mapping of $M$ to a point, then by functoriality we have
$\pi^{\pa M}_*\circ \Theta =\pi^M_*\circ \Theta \circ \iota_*$ as maps
acting on $KK_\Gamma^1(C_0(\pa \wm), C(T)). $
Since $
\iota_*\circ \delta_\Gamma=0$, the previous
lemma implies that
$$
\ind \D_0= \pi^{\pa M}_*\circ \Theta [\D_0]=\pi^{\pa M}_*\circ \Theta \circ
\delta_\Gamma[\D]= \pi^M_*\circ \Theta \circ
\iota_*\circ \delta_\Gamma[\D]=0.
$$ The theorem is proved.
\end{proof}

\subsection{Dirac $b$-index classes on foliated bundles with boundary}\label{subsect:b}

Now let $\P$ be a spectral section for $\D_0$ and
consider an associated trivializing operator $\A^0_\P$
 as in Proposition \ref{Triv}.
Let $\rho \in C^\infty_c(\RR, \RR^+)$ be a nonnegative even smooth
test function such that $\int_{\RR} \rho(x) d x =1.$ We set
$\rho_\epsilon(x) = \frac{1}{\epsilon}\rho(\frac{x}{\epsilon})$
and then consider the Fourier transform of $\rho_\epsilon$: $$
\widehat{\rho_\epsilon}(z) = \int_\RR e^{- i t z}
\rho(\frac{t}{\epsilon} ) \epsilon^{-1} dt. $$ Then there exists a
self-adjoint operator $\A_\P \in \Psi^{-\infty}_{b, \Cr}$ (as already remarked this
space is defined in \cite{LPETALE}) such that the indicial family
of $^{b}\D^+ + \A_\P$ is given by $$ \forall z \in \RR,\;
I(^{b}\D^+ + \A_\P , z) = \D_0 + i z + \widehat{\rho_\epsilon}(z)
\A^0_\P $$ and is invertible from
$H^1_{\Cr}(\pa \wm, \we^+_{| \pa \wm})$ onto
$L^2_{\Cr}(\pa \wm, \we^+_{| \pa \wm})$
 for
any $z \in \RR.$ Recall that $\A_\P$ is constructed in the
following way (see formula (8.7) of \cite{MP I}). Choose a
$\Gamma-$invariant fiberwise product decomposition near the
boundaries of the fibers. Let $\A_\P^\prime$ be the unique
$\RR^+-$invariant operator (Melrose \cite{Melrose} page 126)
 such that:
$$
\forall z \in \RR,\; I( \A_\P^\prime,z)= \widehat{\rho_\epsilon}(z) \A^0_\P.
$$ Then set $\A_\P= \phi(x) \A_\P^\prime \phi(x^\prime)$ where
$\phi\in \C^\infty([0,1], \RR)$ is such that
$\phi(x)=1$ for $x\in [0, \frac{1}{2}]$ and
$\phi(x)=0$ for $x\geq \frac{3}{4}.$

The following theorem is proved
exactly as in Section 3.4 of \cite{LPETALE}.
\begin{theorem} Let $m\in \NN^*.$ The operator $^{b}\D^+ + \A_\P$
defines a $\Cr-$Fred\-holm operator $$H^m_{b , \Cr}( \wm , \we^+)
\longrightarrow H^{m-1}_{b , \Cr}( \wm , \we^-).$$ Its associated index class
does not depend on $m$ and we
 denote it by $\Ind_b
(^{b}\D^+, \P) \in K_0(\Cr).$
\end{theorem}

\medskip
\subsection{Dirac $APS-$index classes   on foliated bundles with boundary}\label{subsect:APS}

We keep the same geometric data as in Subsection
\ref{subsect:preliminaries} but we replace the $\Gamma-$invariant
vertical $b-$metric by a vertical metric having a product
structure near the boundary. We then get a $\Cr-$linear
$\ZZ_2-$graded Dirac type operator $\D$ acting on $C^\infty_c(\wm
, \we):$ $$ \D=\begin{pmatrix}  0 & \D^- \cr \D^+ & 0 \cr
\end{pmatrix}.
$$ One has $\D^+= \sigma ( \frac{\partial}{\partial x} + \D_0 )$ and
$\D^-= \sigma^{-1} ( \frac{\partial}{\partial x} + \sigma \D_0 \sigma^{-1}).$

Consider a spectral section $\P$ for $\D_0$ and define an odd operator
acting on $L^2_{\Cr}(\partial \wm , \we_{\partial \wm})$ by
$$
B_\P= \begin{pmatrix} 0 & (\id -\P) \sigma^{-1} \cr
\sigma^{-1} \P & 0 \cr
\end{pmatrix}
$$ Next we introduce the domain dom$\, \D_{ \P}$ of
$ \D$ associated with the global APS boundary condition defined by $\P:$
$$
{\rm dom}\, (\D_{ \P}) =
\{ \xi \in H^1_{\Cr}(\wm , \we)/\; B_\P(\xi_{|\partial \wm})=0\}
$$ and will denote by $\D_{ \P}$ the restriction of $\D_{}$
to ${\rm dom}\, (\D_{ \P}).$ In a similar and obvious way one defines
$\D_{ \P}^\pm.$

\begin{theorem} \label{Fred} {\item 1)} The operator $\D_{ \P}^+$ defines a
$\Cr-$Fredholm operator from ${\rm dom}\, (\D_{ \P}^+)$
to $L^2_{\Cr}( \wm, \we).$ We denote by
$\Ind^{APS} (\D^+,  \P )\in K_0(\Cr)$ the associated index class.
{\item 2)} One has
$\Ind_b (^{b}\D^+, \P) = \Ind^{APS} (\D^+ , \P ) \in K_0(\Cr).$
\end{theorem}
\begin{proof} 1) The arguments of Wu \cite{Wu I} page 374 can be immediately
extended to our setting and allow to get easily the result.

2) One proceeds exactly as in the proof of Theorem 5 of \cite{LPCUT}.
\end{proof}

\subsection{APS-index theory via the $b$-calculus: an overview of the literature.}\label{subsect:overview}

The use of the $b$-calculus on manifolds with boundary has
generated a great number of interesting articles: in this
subsection we shall review {\it only} those papers that use such a
pseudodifferential calculus and are {\it directly} connected to an
index theorem on manifolds with boundary.

The basic reference for the $b$-calculus on compact manifolds with
boundary is of course the book by Melrose \cite{Melrose}. For a
short introduction to the $b$-calculus and its use in establishing
the APS-index formula the reader can also refer to the surveys of
Mazzeo and Piazza \cite{piazza2} and Grieser \cite{grieser}. The
contribution of Loya in these proceedings \cite{loya2} is also an
excellent introduction. For pseudodifferential extensions of the
Atiyah-Patodi-Singer formula in the context of the $b$-calculus,
one can consult the work of Piazza \cite{piazza-jfa} and
Melrose-Nistor \cite{Melrose-Nistor}. The work of
Nistor-Weinstein-Xu \cite{NWX}  and Monthubert \cite{Monthubert}
in the context of groupoids should also be mentioned. Further
generalizations of the index formula via the $b$-calculus have
been given by  Hassel-Mazzeo-Melrose \cite{asatet} \cite{stmwc} to
manifolds with corners . For more in this direction, see  also the
recent survey article  of Loya \cite{loya}.

The notion of spectral section and its use in establishing a
general APS-{\it family} index theorem for Dirac operators appears
for the first time in the work of Melrose and Piazza \cite{MP I}
\cite{MP II}. Their theorem extends a result of Bismut-Cheeger
\cite{BC2} from the case where the boundary family is invertible
to the general case. For a quick introduction to the results
proved in these articles the reader can  refer to the survey of
Piazza \cite{piazza1}. A pseudodifferential extension of the
Melrose-Piazza family-index theorem has been recently established
by Melrose and Rochon in \cite{Melrose-Rochon}. The family
APS-index theory developed in \cite{MP I} \cite{MP II}
 was extended by
Leichtnam and Piazza to the specific  noncommutative context of
{\it Galois $\Gamma$-coverings with boundary}, see
\cite{LPMEMOIRS} \cite{LPGAFA} \cite{LPCUT}, following a
conjecture of Lott \cite{Lott II}. In these papers not only
suitable index classes are defined in  $K_* (C^*_r \Gamma)$ but
explicit formulae are also obtained for the pairing of these index
classes with suitable cyclic cocycles. {\it Geometric
applications} of these results on Galois $\Gamma$-coverings have
been given to the problem of defining higher signatures on
manifolds with boundary (see Lott \cite{Lott II} \cite{Lott 3})
and proving their homotopy invariance \cite{LLP},  to uniqueness
problems in positive scalar curvature metrics, see Leichtnam-Piazza
\cite{LPPSC}, to the problem of cut-and-paste invariance of
Novikov higher signatures on closed manifolds \cite{LLP},
\cite{LPCUT} (see also \cite{LLK}, Hilsum \cite{Hilsum}), to the
homotopy invariance of the Atiyah-Patodi-Singer and Cheeger-Gromov
rho-invariants for closed compact manifolds having a  torsion-free
fundamental group $\Gamma$ satisfying the bijectivity of the
Baum-Connes map for $C^*_{{\rm max}} \Gamma$ (see Piazza-Schick
\cite{PS}), a result due originally to Keswani \cite{Kes}. The
geometric applications to higher signatures, as well as the index
theorems underlying them, are now also treated in the survey by
Leichtnam-Piazza \cite{LPBOUTET}.

As already explained, in Leichtnam-Piazza  \cite{LPETALE} we
define an index class on foliated bundles with boundary under an
invertibility assumption on the boundary operator; we also
establish an index formula for the higher indeces obtained by
pairing this index class with suitable cyclic cocycles. Finally,
in the present paper we have just defined index classes associated
to an arbitrary Dirac operator on a foliated bundle with boundary
{\it and} the choice of a spectral section for its boundary
operator. We shall now see what are the fundamental properties of
this index class and how they can be employed in order to
investigate the cut-and-paste invariance of the Baum-Connes higher
signatures on closed foliated manifolds.

\section{Fundamental properties of $b$-index classes}\label{sect:fundamentals}

\subsection{The relative index theorem.}

The following theorem extends the special fibration case treated by Melrose and Piazza in \cite{MP I}, as 
well as the covering case in \cite{LPMEMOIRS}. 
\begin{theorem} \label{Rel} Let $\P_1$ and $\P_2$ be two spectral sections
for $\D_0$. Then one has:
$$
\Ind_b (^{b}\D^+, \P_2)-\Ind_b (^{b}\D^+, \P_1)= [\P_1-\P_2] \in K_0(\Cr).
$$
\end{theorem}
\begin{proof}
We shall make precise and at the same time extend the proof sketched 
in \cite{LPAGAG} for $T=$point.

\n
Using Lemma \ref{Fan} and the proofs of Lemma 8 and Proposition 17
of \cite{MP I} one
 checks easily the following five facts:

\n
\item{(a)} One can assume that $\P_1=\id$ on the
range of $\P_2$.

\n
\item{(b)} There exist  two spectral sections
$\Q,\,\R$ for $\D_0$ such that for any $j\in \{1,2\}$:
 $$ \P_j
\Q=\Q \P_j=\Q,\; \P_j \R=\R \P_j= \P_j, \; \Q \R= \R \Q =\Q.
 $$

 \n
\item{(c)} The four following self-adjoint projections
$$\Q\,,\quad \id-\R\,,\quad \P_j\R(\id-\Q)\,,\quad
(\id-\P_j)\R(\id-\Q)
,\; 1 \leq j \leq 2 $$ commute with each
 other and the sum of their (four) ranges provide an orthogonal decomposition
 of $L^2_{\Cr}(\partial \wm ; \we^+_{|\partial \wm})$.

\item{(d)} One has $(\P_1-\P_2) \R (\id -Q)=(\P_1-\P_2).$

\item{(e)} There exists $s>0$ such that for each $j\in \{1,2\}$ the
 operator $\D_0^j :=$
 $$
 \Q \D_0 \Q + s \P_j \R (\id - \Q) +
 (\id-\R)\D_0(\id -\R) -s (\id-\P_j)\R (\id -Q) =
 \D_0 + {\mathcal{A}^0_{\P_j}}
 $$ is invertible.

\m
 Now we set for $r\in [-1,1]$:
 $$
 \D_0(r) = \frac{1}{2} (1+r)\D_0^1 + \frac{1}{2} (1-r)\D_0^2.
 $$
 \begin{lemma} \label{Lem} {\item 1)} For any $r\in [-1,1]\setminus \{0\}$,
 $\D_0(r)$ is invertible.
 {\item 2)} There exists $\epsilon_1> 0$ such that the
  $L^2_{\Cr}-$spectrum of $\D_0(0)$ does not meet
   $]-2 \epsilon_1, 2 \epsilon_1[\setminus \{0\}$,
  $\ker \D_0(0)=[\P_1-\P_2]$ is a $\Cr-$finitely generated
 projective module,
 $$
 L^2_{\Cr}(\partial \wm ; \we^+_{|\partial \wm})= \ker \D_0(0) \oplus
 (\ker \D_0(0))^\perp
 $$ and $\D_0(0)$ defines an $L^2_{\Cr}-$invertible operator from
 $(\ker \D_0(0))^\perp$ onto
 
 \noindent 
 $L^2_{\Cr}(\partial \wm ; \we^+_{|\partial \wm})$.
 \end{lemma}
 \begin{proof} 1) An easy computation shows that
 \begin{equation}\label{Equ}
 \D_0(r)= 
 \Q \D_0 \Q + s \P_2 \R (\id - \Q) +
 (\id-\R)\D_0(\id -\R)\end{equation}
 $$ -s (\id-\P_2-(1+r)(\P_1-\P_2) )\R (\id -Q).$$
 
  Using properties (b) and (c) one gets the result of 1).

 2) The previous identity with $r=0$ and properties (a), (b), (c)
 show that $\ker \D_0(0)$ coincides with the range of
 $(\P_1-\P_2) \R (\id -Q)=(\P_1-\P_2)$ and that
 $\D_0(0)$ defines an invertible operator from
 $(\ker \D_0(0))^\perp$ onto
 $L^2_{\Cr}(\partial \wm ; \we^+_{|\partial \wm})$. One
 then gets immediately part 2).
 \end{proof}

 For each
 $r\in[-1,1]$ set $\D_0(r)= \D_0 + \A^0(r).$
 Just as before Theorem \ref{Fred}, we consider  a
 $\Cr-$linear operator $^{b}\D(r)$ such that:
 $$
 \forall z\in \RR, \quad I(^{b}\D(r)^+,z) = \D_0 + iz \id +\widehat{\rho}_\epsilon(z)
 \A^0(r).
 $$
 Using Lemma \ref{Lem}. 2) and the extension of the Melrose $b$-calculus
 to foliated bundles, as explained in \cite{LPETALE}, one checks easily that for
 $t\in ]0 , \epsilon_1[$,
 $x^{\mp t}\, ^{b}\D(0)^+ x^{\pm  t}$ induces a $\Cr-$Fredholm operator from
 $H^1_{b,\Cr}(\wm, \we^+)$ into $L^2_{b,\Cr}(\wm, \we^-).$
 Denote by $\Ind_{\pm t} ^{b}\D(0)^+$ the corresponding index class
 in $K_0(\Cr).$
 \begin{lemma} \label{Lem1}

 {\item 1)} For any $t\in ]0, \epsilon_1 [$ one has:
 $$\Ind_{- t} \,^{b}\D(0)^+=\Ind\,^{b}\D(-1)^+.
 $$
 { \item 2)}  For any $t\in ]0, \epsilon_1 [$ one has:
 $$
 \Ind_{ t} \,^{b}\D(0)^+=\Ind \, ^{b}\D(1)^+.
 $$
 \end{lemma}
 \begin{proof} 1) Fix $ t\in ]0, \epsilon_1[.$
 For any $r\in [-1,0]$ one has:
 $$
 \forall z\in \RR, \quad I( x^{t}\, ^{b}\D(r)^+ x^{-t},z)  = \D_0 + (iz -t)\id +\widehat{\rho}_\epsilon(z)
 \A^0(r).
 $$ Observe that for any $r\in [-1,0]$, $-t-s +s(1+r)\not=0.$
 Using Lemma \ref{Lem} and inspecting expression  \eqref{Equ}  (especially
 its last term), one checks
 immediately that  for any $r\in [-1,0]$, $I( x^{t}\, ^{b}\D(r)^+ x^{-t},0)$
 is $L^2_{\Cr}-$invertible.
 Recall
 that the $\widehat{\rho}_\epsilon(z)$ take real values, then since
 $\D_0$ and $\A^0(r)$ are self-adjoint it is clear
 that for any $(r, z)\in [-1,0]\times \RR^*$, $I( x^{t}\, ^{b}\D(r)^+ x^{-t},z)$
 is $L^2_{\Cr}-$invertible.Therefore,
 for any $t \in ]0, \epsilon_1 [$ the family
 $\{x^{t}\, ^{b}\D(r)^+ x^{-t},\; r\in [-1, 0]\}$ defines a continuous family
 of $\Cr-$Fredholm operators. By the homotopy invariance
 of the $\Cr-$index in $K_0(\Cr)$ one has:
 $$
 \forall r\in [-1, 0],\;\; \Ind\,  x^{t}\, ^{b}\D(r)^+ x^{-t}\,=\,
 \Ind\,  x^{t}\, ^{b}\D(0)^+ x^{-t}.
 $$ Next, the family $\{x^{t^\prime}\, ^{b}\D(-1)^+ x^{-t^\prime},
 \; t^\prime\in [0, t]\}$ defines a continuous family
 of $\Cr-$Fredholm
 operators so that one gets:
 $$
  \Ind\,x^t\,^{b}\D(-1)^+x^{-t} \,=\,\Ind\,^{b}\D(-1)^+.
 $$ From the last two equations one gets immediately part 1).
 Part 2) of the lemma is
  proved in a similar way.
 \end{proof}
 Since one has:
 $$
 \Ind\, (^{b}\D^+, \P_2)= \Ind \,^{b}\D(-1)^+,\;\,
 \Ind \, (^{b}\D^+, \P_1)= \Ind\,  ^{b}\D(1)^+
 $$ Lemma \ref{Lem1} shows that Theorem \ref{Rel} is a consequence
 of the following proposition
 \begin{proposition} For any $t\in ]0, \epsilon_1 [$, one has:
 $$\Ind_{- t} \,^{b}\D(0)^+-\Ind_{ t} \,^{b}\D(0)^+=[\ker \D_0(0)]=
 [\P_1-\P_2].
 $$
 \end{proposition}
 \begin{proof} Denote by $\D_{cyl}^+$ the $\Gamma-$equivariant family
 of fiberwise elliptic operators which in a collar neighborhood ($\simeq
 [0,1]\times \pi^{-1}(\theta)$) of the
 boundaries is given by
  $$\D_{cyl}^+= \sigma ( x \partial_x + \D_0(0))$$ and which
  coincides with $^{b}\D^+$ outside this collar neighborhood.
  Observe that unlike $^{b}\D^+$, $\D_{cyl}^+$ is not a $b-$operator.
  Proceeding as in Section 10 of \cite{LLP}, one proves that
  $\D_{cyl}^+$ acting from $x^{\pm t}H^1_{b,\Cr}(\wm, \we^+)$
  into $x^{\pm t}L^2_{b,\Cr}(\wm, \we^-)$ defines a $\Cr-$Fredholm operator
  whose index class $\Ind_{\pm t} \, \D_{cyl}^+$ satisfies:
  $\Ind_{\pm t} \, \D_{cyl}^+= \Ind_{\m t} \,^{b}\D(0)^+.$

  Now, proceeding as in the proof of Proposition 4 of
 \cite{MP I} one checks that there exists a positive number
 $N$ and two continuous $\Cr-$linear maps
 $\R_{\pm t}: (\Cr)^N \rightarrow x^{2t} H^{2n+3}_{b,\Cr}(\wm, \we^-)$ such that
 the following two maps are surjective:
 $$
 \D^t:  x^t H^1_{b,\Cr}(\wm, \we^+)\oplus (\Cr)^N \rightarrow
 x^t L^2_{b,\Cr}(\wm, \we^-)
 $$
 $$
 \D^{-t}:  x^{-t} H^1_{b,\Cr}(\wm, \we^+)\oplus (\Cr)^N
 \rightarrow
 x^{-t} L^2_{b,\Cr}(\wm, \we^-).
 $$ 
 with
 $$ \D^t:= \D_{cyl}^+ + \R_t\,,\quad \D^{-t}= \D_{cyl}^+ + \R_{-t}\,.$$
 Then one has:
 $$
 \Ind_{ t} \D_{cyl}^+= [\ker \D^t] - [(\Cr)^N],\;
 \Ind_{ -t} \D_{cyl}^+= [\ker \D^{-t}] - [(\Cr)^N].
 $$ Thus we just have to prove that:
 $$
 [\ker \D^{-t}] - [\ker \D^t]= [\P_1-\P_2]=\ker \, \D_0(0) \in K_0(\Cr).
 $$ We are going now to show the existence of a short
 exact sequence
 $$0\rightarrow \ker \D^{ t} \rightarrow \ker \D^{- t} \rightarrow \ker \D_0(0)
  \rightarrow 0.
 $$ Since $\ker \D_0(0)$ is projective we shall obtain that
 $\ker \D^{- t}= \ker \D^{ t} \oplus \ker \D_0(0)$ which will
 imply the proposition.
Our arguments are very much inspired by those used by Melrose in
his proof of the relative index formula in \cite{Melrose}.

 Consider $(u\oplus a ) \in  x^{-t} H^1_{b,\Cr}(\wm, \we^+)\oplus (\Cr)^N$
 such that $\D^{-t}(u\oplus a)=0.$ Let $\phi\in C^\infty([0,1], \RR)$
 be such that $\phi(x)=1$ for $x\leq \frac{1}{2}$ and
 $\phi(x)=0$ for $x\geq \frac{3}{4}.$ Then
 $\D^{-t}(\phi \cdot u) \in x^{2 t}L^2_{b, \Cr}(\wm, \we^-).$
 Denote by $U_M(z,y)$ the Mellin transform of $\phi \cdot u:$
 $$
 \forall y \in \pa \wm,\; U_M(z,y)=\int_{\RR} x^{-i z} (\phi \cdot u)(x,y) \frac{d x}{x}.
 $$ We observe that $z\rightarrow U_M(z,y)$
 is holomorphic
 on the half plane $\{\Im z > t\}$.
 Moreover, since
 $$
 \D^{-t}(\phi \cdot u) = \sigma ( x \frac{ \partial}{\partial x} +
 \D_0(0) ) (\phi \cdot u)\in x^{2 t}L^2_{b, \Cr}(\wm, \we^-)
 $$ one checks easily that
 $z\rightarrow (i z + \D_0(0)) U_M(z,y)$
 is holomorphic
 on the half plane $\{\Im z > -2 t \} .$
  We recall the orthogonal decomposition of Lemma \ref{Lem}   2), and  write for
  each $z \in \{ z^\prime,\; \Im z^\prime > - 2 t \}$:
  $$
 (i z + \D_0(0)) U_M(z,y) =W_0(z,y) \oplus W_1(z,y)
  \in \ker \D_0(0) \oplus (\ker \D_0(0))^\perp.
  $$ Then for any $z \in \{ z^\prime,\; \Im z^\prime > t \}$:
  $$
  U_M(z,y) = \frac{1}{i z} W_0(z,y) \oplus (i z + \D_0(0))^{-1}W_1(z,y).
  $$

Considering the inverse Mellin transform one checks easily
  that
  $$
  \frac {\phi(x)}{ 2 \pi} \int_{\Im z = - \frac{3 t}{2}}
  x^{ i z} \left( \frac{ W_0(z,y) - e^{-z^2}W_0(0,y)}{i z} +
  (i z + \D_0(0))^{-1}W_1(z,y) \right) dz \,
  $$ belongs to $ x^t H^1_{b,\Cr}(\wm, \we^+) .$ Similarly, one checks
  that
  $$
  \frac {\phi(x)}{ 2 \pi} \int_{\Im z = - \frac{3 t}{2}}
  x^{ i z} ( \frac{ W_0(0,y) }{i z}
   ) dz \, \in x^t H^1_{b,\Cr}(\wm, \we^+)
  $$ if and only if $W_0(0,y)\equiv 0.$
   Then set:
  $$
  \Pi_0(u) = W_0(0,y).
  $$
  From our previous computations one gets the
   following lemma.
  \begin{lemma} \label{LEM} With the previous notations:
  $u\in x^t H^1_{b,\Cr}(\wm, \we^+)$ if and only if
  $\Pi_0(u) = W_0(0,y)$  is the null element of $\ker \D_0(0).$ Moreover,
  the following sequence
  $$
  0\rightarrow \ker \D^{ t} \rightarrow \ker \D^{- t} \overset{\Pi_0}
  \rightarrow \ker \D_0(0)
  $$ is exact.
  \end{lemma}

  Now we are going to show that the map $ \ker \D^{- t} \overset{\Pi_0}
  \rightarrow \ker \D_0(0)$ is surjective.

  Consider an element $V_0(y)$ of $\ker \D_0(0)$ and  set:
  $$
  v_1(x,y) = \frac {\phi(x)}{ 2 \pi} \int_{\Im z =  t}
  x^{ i z} \frac{ e^{-z^2} }{ i z} V_0(y) d z.
  $$ It is clear that
  $$
  \sigma^+ ( x \pa_x + \D_0(0)) v_1(x,y) \in x^{2 t} L^2_{b, \Cr}( \wm, \we^-),
  $$ then since $\D^t$ is surjective, there exists
  $u_2 \oplus a_2 \in x^{t} H^1_{b,\Cr}(\wm, \we^+)\oplus (\Cr)^N$ such that
  $\D^t(u_2 \oplus a_2) = - \D^{-t}(v_1 \oplus 0).$ Then
  $(v_1 -u_2) \oplus (-a_2) $ belongs to
  $\ker \D^{- t}$ and $\Pi_0(v_1-u_2)  = V_0(y)$.
  Then, using Lemma \ref{LEM} one obtains (as previously announced) the following
  short exact sequence of $\Cr-$finitely generated projective modules:
  $$
  0\rightarrow \ker \D^{ t} \rightarrow \ker \D^{- t} \overset{\Pi_0}
  \rightarrow \ker \D_0(0)
  \rightarrow 0.$$
  Since the modules are projective one obtains an isomorphism:
$\ker \D^{- t} \simeq \ker \D^{ t} \oplus \ker \D_0(0)$ from which
the proposition follows.
  \end{proof}

 Theorem \ref{Rel} is thus proved.
\end{proof}

\medskip
\subsection{The gluing formula.}

 Let $T$ be
a smooth closed compact connected manifold on which $\Gamma$ acts
on the right. Let $\widehat{M}$ be a closed manifold  on which
$\Gamma$ acts freely, properly and cocompactly on the right: the
quotient space $M=\widehat{M} /\Gamma$ is thus a smooth compact
manifold. We assume that $\widehat{M}$ fibers over $T$ and that
the resulting fibration $$ \pi: \,\widehat{M} \rightarrow T $$ is
a {\it $\Gamma$-equivariant fibration} with fibers
$\pi^{-1}(\theta),\theta \in T,$  of dimension $2k$. Notice that
each fiber is  smooth; we shall also denote the typical fiber of
$\pi: \,\widehat{M} \rightarrow T $ by $Z$. We choose a
$\Gamma$-invariant metric
 on the vertical tangent
bundle $ TZ $. Finally, we  assume
the existence of a
$\Gamma-$equivariant spin structure on $ TZ $
that is fixed once and for all.
We denote  by
$S^Z \rightarrow \widehat{M} $ the associated spinor bundle.

\smallskip
We consider also a $\Gamma$-equivariant complex hermitian vector bundle
$\widehat{V} \rightarrow \widehat{M} $ endowed with a
$\Gamma-$invariant hermitian connection $\widehat{\nabla}.$
We then set $\widehat{E}= S^Z \otimes \widehat{V}= \widehat{E}^+ \oplus
\widehat{E}^-$ which defines a smooth $\Gamma-$invariant
 family of  $\ZZ_2-$graded hermitian
Clifford modules on the fibers $\pi^{-1}(\theta), \, \theta \in T$.
We then get a smooth family of $\Gamma-$invariant  $\ZZ_2-$graded Dirac
type operators
$$  D(\theta)= \begin{pmatrix} 0 & {D^-(\theta)}
\cr {D^+(\theta)}  & 0 \cr
\end{pmatrix},\; \theta\in T
$$
 acting fiberwise on
$C^\infty_c( \widehat{M},\, \widehat{E})$.

\smallskip The family $(D(\theta))_{\theta \in T}$ defines a
$\Cr-$linear $\ZZ_2-$graded Dirac type operator $\D$ acting
on $C^\infty_c( \wm , \we).$ This operator has a well defined
index class in $K_0( C(T) \rtimes_r \Gamma).$ Now let $F$
be a closed cutting $\Gamma-$invariant hypersurface of $M$ such that
$\widehat{M}=\widehat{M}_+ \cup \widehat{M}_-$ where $\widehat{M}_\pm$
are two manifolds whose common boundary is $F$ and which
both fiber over $T$: $\widehat{M}_\pm \rightarrow T.$ We assume that all
these data have a product structure near $F$. Let $\P$ and $\Q$
be two spectral sections for the boundary operator
of the operator induced by the restriction $\D_{|\widehat{M}_+}$ of
$\D$ to $\widehat{M}_+.$ Observe that $\id -\Q$ is a spectral
section for the boundary operator of $\D_{|\widehat{M}_-}$. The following
gluing formula is proved exactly as page 380 in \cite{LPCUT} using an idea
of U. Bunke \cite{Bunke}.
It generalizes the result of Dai-Zhang \cite{xdai} in the fibration case as well as the result 
on covering space in \cite{LPMEMOIRS}. 
\begin{theorem}
$$\Ind \D^+ =
\Ind^{APS} ( \D^+_{|\widehat{M}_+},  \P ) +
\Ind^{APS} ( \D^+_{|\widehat{M}_-}, \id -\Q) + [\P-\Q].
$$
\end{theorem}

\subsection{The variational formula.}

Let $\pi: \widehat{M}\rightarrow T$
be $\Gamma-$equivariant fibration (whose fibers are manifold with boundary)
exactly as in Section 4.1.
We assume that there exists a smooth $1-$parameter family of fiberwise vertical
$\Gamma-$invariant riemannian metrics
$(g_u)_{u \in [1,2]}$ on $\widehat{M}$ which all have a product structure
near the boundary. We assume that the $\ZZ_2-$graded Clifford module
$\widehat{E}$ is endowed with continuous $1-$parameter families
of $\Gamma-$equivariant hermitian metrics $(h^u)_{u \in [1,2]}$
and $\Gamma-$equivariant hermitian connection $(\nabla^u)_{u \in [1,2]}.$
We denote by $(\D_u)_{u \in [1,2]}$ the associated family
of $\Cr-$linear Dirac type operator acting
on $H^1_{\Cr}(\wm , \we).$ We can then state the following
variational formula which is (as in \cite{LPCUT} page 383) an easy
consequence of the relative index theorem Theorem \ref{Rel}.
It generalizes the result of Dai-Zhang \cite{xdai}  in the fibration case as well as the result 
on covering space in \cite{LPMEMOIRS}. 
\begin{proposition} Let us denote by
$\{(\D_u)_0,\; u \in [1,2]\}$ the family of boundary operators
associated to $(\D_u)_{u \in [1,2]}$. We fix noncommutative spectral sections
$\P_1, \P_2$ for $(\D_1)_0$ and $(\D_2)_0$ respectively. Then:
$$
\Ind^{APS} (\D^+_2, \P_2) - \Ind^{APS}(\D^+_1, \P_1)
= {\rm sf}(\{(\D_u)_0\}; \P_2, \P_1) \; {\rm in}\; K_0(\Cr).
$$
\end{proposition}


\section{On the cut-and-paste invariance of the signature index class}
\label{sect:cut-and-paste}

\subsection{Cut-and-paste on foliated bundles.}
\label{subsect:Cut-and-paste}

 We first consider a
$\Gamma-$equivariant fibrations $\pi_{\widehat{X}}: \widehat{X}
\rightarrow T$ with oriented fibers
 and such that the quotient $X=\widehat{X}/\Gamma$
is a smooth compact manifolds.
 Let $r :X=
\widehat{X}/\Gamma \rightarrow (E\Gamma \times T)/\Gamma $ be the
classifying map of the action of the groupoid $T\rtimes \Gamma$ on
$\widehat{X}$ (Connes, \cite[Chapter III]{Co}, \cite{Go-Lo}). This map is
the defined as follows: the $\Gamma$-covering
 $\rho:\widehat{X}\rightarrow X=\widehat{X}/\Gamma $ is
classified by a $\Gamma$-equivariant map $\tilde{\rho}:
\widehat{X}\rightarrow E\Gamma $; let
$\widehat{r}:\widehat{X}\rightarrow E\Gamma\times T$ be the map
$(\tilde{\rho},\pi_{\widehat{X}})$; then $r :X=
\widehat{X}/\Gamma \rightarrow (E\Gamma \times T)/\Gamma$ is the
$\Gamma$-quotient of $\widehat{r}$.

 We now consider
two $\Gamma-$equivariant fibrations $\pi_{\widehat{M}}:
\widehat{M} \rightarrow T$ and $\pi_{ \widehat{N}}: \widehat{N}
\rightarrow T$ where in both cases the fibers are even
$2m-$dimensional oriented manifolds with boundary and such that
the quotient $M=\widehat{M}/\Gamma$ and $N=\widehat{N}/\Gamma$ are
two smooth compact manifolds with boundary.

We assume the existence of two $\Gamma-$equivariant
diffeomorphisms $\phi,\psi:
\partial \widehat{M} \rightarrow \partial \widehat{N}$ such that $
\pi_{\widehat{\partial N}}\circ \phi= \pi_{\widehat{\partial M}}$,
$ \pi_{\widehat{\partial N}}\circ \psi= \pi_{\widehat{\partial
M}}$ and $\phi$, $\psi$ both preserve the orientations of the
fibers. We set:
$$
\widehat{X}_\phi=\widehat{M}
\cup_\phi\widehat{N}^-,\; \widehat{X}_\psi=\widehat{M}
\cup_\psi\widehat{N}^-
$$ where $\widehat{N}^-$ means that the fibers of $\widehat{N}$ are endowed 
with the reverse orientation.

One then  consider the two $\Gamma-$equivariant
fibrations $ \pi_\phi:\widehat{X}_\phi\rightarrow T$ and
$\pi_\psi: \widehat{X}_\psi \rightarrow T$ and the
two $\Gamma-$covering maps:
$$
\rho_\phi:\widehat{X}_\phi \rightarrow \widehat{X}_\phi/\Gamma,\;
\rho_\psi:\widehat{X}_\psi \rightarrow \widehat{X}_\psi/\Gamma.
$$ Denote by $\tilde{\rho}_\phi: \widehat{X}_\phi \rightarrow E \Gamma$
(resp. $ \tilde{\rho}_\psi: \widehat{X}_\psi \rightarrow E \Gamma$)
the corresponding classifying map of $\rho_\phi$  (resp. $\rho_\psi$).
Then, as explained above, the pairs $( \tilde{\rho}_\phi,\pi_\phi )$ and
$(\tilde{\rho}_\psi ,\pi_\psi )$ induce two
classifying maps:
$$ r_\phi :\widehat{X}_\phi\,/\Gamma \rightarrow (E\Gamma \times
T)/\Gamma,
 \; s_\psi :\,\widehat{X}_\psi/\Gamma \rightarrow (E\Gamma \times T)/\Gamma. $$

We shall briefly say that two
 $T\rtimes \Gamma$-proper manifolds obtained  above
  are {\it cut-and-paste
equivalent}.


\subsection{The defect formula.}\label{subsect:defect}

  We consider
a fiberwise vertical metric $g_1$ (resp. $g_2$) on
$\widehat{X}_\phi$ (resp. $\widehat{X}_\psi$) which is product
like near $\partial \widehat{M}.$ Consider $(g_1)_{|\widehat{M}}$
and $(g_2)_{|\widehat{M}}$ and let $g_{+,u},$ with $u \in [1,2],$
be a path of vertical fiberwise riemannian metrics on
$\widehat{M}$ connecting them and having a product structure near
the boundary. Similarly we choose a path $(g_{-,u})_{u\in [1,2]}$
of vertical fiberwise riemannian metrics on $\widehat{N}$
connecting $(g_1)_{|\widehat{N}}$ and $(g_2)_{|\widehat{N}}.$ One
thus gets two family of boundary $\Cr-$linear signature  operators
(as in Section 4.2) $\{\D^{{\rm sign} ,u}_{\partial \widehat{M}}\}_{u\in [1,2]}$
and $\{\D^{{\rm sign} ,u}_{\partial \widehat{N}}\}_{u\in [1,2]}.$ Observe that
$\D^{{\rm sign} ,1}_{\partial \widehat{M}}$ is conjugated through $\phi^*$ to
$-\D^{{\rm sign} ,2}_{\partial \widehat{N}}$ and that $\D^{{\rm sign} ,2}_{\partial
\widehat{M}}$ is conjugated through $\psi^*$ to $-\D^{{\rm sign} ,2}_{\partial
\widehat{N}}$. We can thus (as in Section 6.1 of \cite{LPCUT}) put
together the family $\{\D^{{\rm sign}, u}_{\partial \widehat{M}}\}_{u\in [1,2]}$
and the family $\{\D^{{\rm sign}, u}_{\partial \widehat{N}}\}_{u\in [1,2]}$ and
obtain a family $$ \{ \D^{{\rm sign}}_{\partial \widehat{M}}(\theta) \}_{\theta
\in S^1}= \{\D^{{\rm sign}, u}_{\partial \widehat{M}}\}_{u\in [1,2]} \cup
\{\D^{{\rm sign}, u}_{\partial \widehat{N}}\}_{u\in [2,1]} $$ which is an
$S^1-$family acting on the fibers of the mapping torus 
defined by $\phi^{-1} \circ \psi.$ Then one has the following
formula  whose proof is an easy extension of the one of Theorem 11
of \cite{LPCUT} which is in turn modeled on the arguments given in
Section \ref{sect:sflow}.

\begin{theorem}\label{defect}  Denote by $\D^{{\rm sign}}_{\widehat{X}_\phi}$ 
(resp. $\D^{{\rm sign}}_{\widehat{X}_\psi}$)
the $\Cr-$linear signature  operator of $\widehat{X}_\phi$ (resp.
$\widehat{X}_\psi$) defined as in Section {\rm 4.4}. The following
formula holds $$ \Ind \D^{{\rm sign}}_{\widehat{X}_\phi} - \Ind
\D^{{\rm sign}}_{\widehat{X}_\psi} = {\rm sf} (\{ \D^{{\rm sign}}_{\partial
\widehat{M}}(\theta) \}_{\theta \in S^1} )\quad \text{in} \quad
K_0(\Cr)\,. $$
\end{theorem}

\subsection{Vanishing spectral flow.}\label{subsect:vanishing-higher}

Let the fibers of $Z\rightarrow \widehat{M} \rightarrow T$
have dimension $2m$, so that the fibers of the boundary fibration
$\pa\widehat{M}\to T$ have dimension $2m-1$. 
We endow the boundaries of the fibers of
$\pi_{\widehat{M}}:  \widehat{M} \rightarrow T$ with a
$\Gamma-$invariant metric and make the following "middle-degree"
assumption on the boundary:
\begin{ass} \label{Lottgroupoid} There exists $\epsilon \in ]0,1[$ such that for each
$\theta \in T$, the $L^2-$spec\-trum of the fiberwise differential-form
laplacian  acting on $$L^2(\partial  \pi_{\widehat{M}}^{-1}(\theta) ;
\wedge^{m-1}T^*\partial  \pi_{\widehat{M}}^{-1}(\theta) )$$
does
not meet $]-\epsilon, \epsilon[$.
\end{ass}

We give an example (inspired by
\cite{LLP} page 563) where this assumption is satisfied.
Let $\widetilde{N} \rightarrow N$ a Galois $\Gamma-$covering
of a smooth orientable compact $2m-$dimensional manifold
with boundary such that $\pa N$ has a cellular decomposition
without any cells of dimension $m$. Set
$\widehat{M}=\widetilde{N} \times T$ consider the
trivial fibration
$$
\pi_{\m}: \widetilde{N} \times T \rightarrow T,\; \pi(z,\theta)=\theta.
$$ Then the above assumption is satisfied in this case.

\n
{\bf Remarks.} 

\noindent 1) In the Galois-covering case ($T=$point), this assumption comes
from the work of John Lott, see \cite{Lott II}.

\noindent  2) Under the
Assumption \ref{Lottgroupoid} one can prove easily that the index
of the boundary $\Cr-$linear signature operator vanishes in
$K_1(\Cr)$.

\noindent 3) Proposition 1 of \cite{LPCUT} shows that
the associated $\Cr-$linear  signature-laplacian  
$(\D^{{\rm sign}}_{\pa \wm})^2$ 
 induces an invertible operator
from the Hilbert $\Cr$-module $H^2_{\Cr}( \pa \wm , \wedge^{m-1} T^* \pa \wm)$ onto $L^2_{\Cr}(
\pa \wm , \wedge^{m-1} T^* \pa \wm).$

\begin{proposition} \label{Prop}
Let $\widehat{X}_\phi$ and $\widehat{X}_\psi$ be as in Subsection
{\rm \ref{subsect:Cut-and-paste}}. 
Denote by $\D^{\rm sign}_{\widehat{X}_\phi}$ (resp.
$\D^{\rm sign}_{\widehat{X}_\psi}$) the corresponding $\Cr-$linear signature operator
of $\widehat{X}_\phi$ (resp. $\widehat{X}_\psi$). Assume that
Assumption {\rm \ref{Lottgroupoid}} is satisfied for $\partial
\widehat{X}_\phi$ 
instead of $\partial \widehat{M}.$ Then
one has:
 $$\Ind \D^{\rm sign}_{\widehat{X}_\phi}
= \Ind \D^{\rm sign}_{\widehat{X}_\psi} \in K_0 ( \Cr)\otimes_\ZZ \QQ.$$
\end{proposition}
\begin{proof} We follow pages 391-392 of \cite{LPCUT}.
We denote by $Z$ the typical fiber of the fibration $\pi: \partial
\widehat{M} \rightarrow T$ and set $\Omega^\ast =
\cap_{j\geq 0}  H^j_{\Cr}(\partial \wm ; \wedge^\ast T^* Z).$
 We then set: $$ V=
d^* \Omega^m + d \Omega^{m-1},\; W=\Omega^< \oplus \Omega^>,\;
{\rm where} $$ $$ \Omega^<= \Omega^0 \oplus \ldots \oplus
\Omega^{m-2} \oplus (d^*\Omega^m)^\perp,\; \Omega^>= (d
\Omega^{m-1})^\perp \oplus \Omega^{m+1}\oplus \ldots \oplus
\Omega^{2m-1}. $$ It is clear that the $\Cr-$linear signature
operator $\D^{\rm sign}_{\partial \wm}$ of $\partial \wm$ sends $V$ (resp.
$W$) into itself. Using Assumption \ref{Lottgroupoid} one checks
easily that $\D^{\rm sign}_{\partial \wm}$ induces an invertible operator on
the $L^2_{\Cr}-$completion of $V$ with domain $H^1$ and we denote
by $\Pi_>$ the projection onto the positive part. Then, proceeding
as in page 392 of \cite{LPCUT}, one checks that $\D^{\rm sign}_{\partial
\wm}$ admits a symmetric spectral section $\P$ in the sense that
$\P$ is diagonal with respect to the decomposition
$\sum_{j=0}^{2m-1} \Omega^j = V\oplus W$ and $$ \P_{|V}= \Pi_>,\;
\, \alpha \circ \P_{|W} + \P_{|W} \circ \alpha =\alpha $$ where
$\alpha$ is the involution of $W$ equal to the identity of
$\Omega^<$ and to minus identity on $\Omega^>.$ Moreover,
proceeding as in the proof of Proposition 17 of \cite{LPAGAG}, one
checks that if $\Q$ is another symmetric spectral section for
$\D^{\rm sign}_{\partial \wm}$ then one has $[\P-\Q]=0$ in
$K_0(\Cr)\otimes_\ZZ \QQ.$ Lastly we observe that in the
definition of the spectral flow ${\rm sf}(\{ \D^{\rm sign}_{\pa
\wm}(\theta)\})_{\theta \in S^1}  $ we may assume that all the
spectral sections involved are symmetric. Then the result follows
from Theorem \ref{defect}.
\end{proof}

\section{Geometric
applications}\label{sect:applications}

In all this section
we shall assume, for simplicity, that $T $ is orientable and that
$\Gamma$ preserves the orientation of $T$.
We shall first define the Baum-Connes higher signatures
of a  $T\rtimes \Gamma-$proper manifold; these are numeric invariants.
 Then we shall ask ourselves
when these higher signatures are cut-and-paste invariant.
The strategy here is to use the basic assumption
\ref{Lottgroupoid} and the equality of the signature index classes
established in Proposition \ref{Prop} in order to deduce
the equality of the higher signatures (or, at least, the equality
of some of these higher signatures).
When $T=$point there are two techniques allowing to use the equality
of index classes in order to deduce the equality of (all) higher signatures:
the first one employs cyclic cohomology and the second one employs
the assemply map from topological $K$-homology to the K-Theory
of the reduced $C^*$-algebra. 
We shall try to generalize these two approaches.

\subsection{Baum-Connes higher signatures.}

 We consider a $(T\rtimes \Gamma)$-proper manifold, i.e. a
$\Gamma-$equivariant fibration  $\pi_{\widehat{X}}: \widehat{X}
\rightarrow T$ with oriented fibers
 and such that the quotient $X=\widehat{X}/\Gamma$
is a smooth compact manifolds.
 Let $r :X=
\widehat{X}/\Gamma \rightarrow (E\Gamma \times T)/\Gamma $ be the
associate classifying map .

For each cohomology class $c\in H^{*}( (E\Gamma \times T)/\Gamma , \CC)$
the number \begin{equation} \label{signature}
 \int_X L(TX) \wedge r^*(c)
 \end{equation} is called the
Baum-Connes higher signature of the $T\rtimes \Gamma-$proper
manifold $ \widehat{X}$ (see Baum-Connes \cite{Baum-Connes}).
We are interested in the set  $$
\{ \int_X L(TX) \wedge r^*(c), \quad c\in H^{*}( (E\Gamma \times T)/\Gamma , \CC)
\}.$$

As already explained, projecting the fibers
$\pi_{\widehat{X}}^{-1} (\theta) $ ($\theta \in T$) onto the quotient
 $X:=\widehat{X}/\Gamma$ one gets
a foliation $\mathcal{F}$ of the compact manifold
$X$. Then the sets of higher signatures
$$
\{ \int_X L(TX) \wedge r^*(c), \quad c\in H^{*}( (E\Gamma \times T)/\Gamma , \CC)
\}$$
 can equally be described as the set
$$\{ \int_X L(T \mathcal{F}) \wedge r^*(c),\quad
[c]\in H^{*}( (E\Gamma \times T)/\Gamma , \CC) \},$$
(here we use the fact that the
$L$-class of the normal bundle to the foliation is
 the pull-back of a class in $
H^{*}( (E\Gamma \times T)/\Gamma , \CC)$).

Next, with the notations of Subsection \ref{subsect:Cut-and-paste},
 let us mention that
 the
 higher signatures of  two cut-and-paste equivalent
$T\rtimes \Gamma-$proper manifolds $\widehat{X}_\phi
$ and  $\widehat{X}_\psi$:
 $$\int_{ {\widehat{X}_\phi }/{\Gamma}} L({\widehat{X}_\phi
 }/{\Gamma} )\wedge r_\phi^*(c),\;\quad
\int_{ {\widehat{X}_\psi }/{\Gamma}} L(
{\widehat{X}_\psi }/{\Gamma}) \wedge
s_\psi^*(c) $$ do {\it not} coincide in general. See Karras-Kreck-Neumann-Ossa
\cite{Karras-Kreck-Neumann-Ossa (1973)} and \cite{LLK} for
examples when $T$ is reduced to a point.

Our goal is to find sufficient conditions on the group $\Gamma$
and on its action, ensuring that the Baum-Connes
higher signatures are cut-and-paste invariant.

\subsection{Cut-and-paste invariance: the cyclic cohomology approach.}

We begin by  recalling several results from Connes \cite{Co} and Gorokhovsky-Lott \cite{Go-Lo}.
Set
$$
{\mathcal{B}}^\omega=\{ \sum_{\gamma \in \Gamma} c_\gamma \gamma:\,
|c_\gamma|\, {\rm decays}\, {\rm faster}\, {\rm than}\, {\rm any}\,
{\rm exponential}\, {\rm in}\, ||\gamma|| \}.
$$
In Section 3 of \cite{Go-Lo} is defined a certain
algebra of noncommutative differential forms
$\Omega^*(T, {\mathcal{B}}^\omega)$.
Consider
a closed graded $N-$trace $\eta$ on $\Omega^*(T, {\mathcal{B}}^\omega)$ which
is concentrated at the identity conjugacy class of $\Gamma$; then
to $\eta$ one associates a cohomology class $ \Phi_\eta \in
H^{N + {\rm dim}\, T + 2 \ZZ}( (E\Gamma \times T)/\Gamma ; \CC)$
where $2 \ZZ$ denotes an even-odd grading.
In fact, one of the main  results of  \cite{Go-Lo} is the
proof of an isomophism between the homology of closed graded $N-$traces
on $\Omega^* (T,\CC\Gamma)$,
concentrated at the identity conjugacy class, and the cohomology space
$H^{N + {\rm dim}\, T + 2 \ZZ}( (E\Gamma \times T)/\Gamma ; \CC)$.
The above class is obtained by restricting our $\eta$
from $\Omega^*(T, {\mathcal{B}}^\omega)$
to
$\Omega^* (T,\CC\Gamma)$ and applying the isomorphism.

\s
{\it From now on all our graded traces will be concentrated
at the identity conjugacy class.}

\s
Let now $\widehat{X}\to T$ a $(T\rtimes \Gamma)$-proper manifold.
Under the hypothesis and notations of Subsection \ref{subsect:Cut-and-paste},
  Gorokhovsky and Lott give a heat-equation proof of Connes' index
theorem for $(T\rtimes \Gamma)$-proper manifolds (see \cite{Co}
Theorem 12, Ch III, Section 7.$\gamma$ for the original
statement). In particular Gorokhovsky and Lott
prove the following formula:
$$
\langle {\rm ch}\, \Ind_\omega \D^{\rm sign} , \eta \rangle = \int_{\widehat{X}/\Gamma}
L(T {\mathcal{F}} ) \wedge r^*(\Phi_\eta).
$$
The index class in this formula is not the index class
defined in Subsection \ref{index-classes-closed}  but rather a refinement
in the $K$-theory of the algebra $$C^\infty (T,\B^\omega)=
\B^\omega\otimes_{\CC\Gamma} (C^\infty_c (T)\rtimes \Gamma),$$
a subalgebra of $C(T)\rtimes_r \Gamma$ containing
$C^\infty_c (T)\rtimes \Gamma$.
 In general, the $K$-theory
groups of $C^\infty (T,\B^\omega)$ and $C(T)\rtimes_r \Gamma$ are different
and the index class defined in 
Subsection \ref{index-classes-closed}  is  the image
of $\Ind_\omega \D^{\rm sign}$ under the $K$-theory homomorphism
$K_0 (C^\infty (T,\B^\omega))\to K_0 (C(T)\rtimes_r \Gamma)$
induced by the inclusion
$C^\infty (T,\B^\omega)\hookrightarrow C(T)\rtimes_r \Gamma$.
On the other hand,
all our formulas have been established for
index classes in  $K_* (C(T)\rtimes_r \Gamma )$;
this means that we need further hypothesis in order to
combine the Gorokhovsky-Lott index formula and
our results in the previous sections.
The result we need is stated in Corollary 3 of \cite{Go-Lo}:
assume that $\mathcal{A}$ is a dense holomophically closed
subalgebra of $C(T)\rtimes_r \Gamma$
containing $C^\infty (T,\B^\omega)$; then we know
that $K_* (\mathcal{A})\simeq K_* (C(T)\rtimes_r \Gamma )$.
Let $\eta$ be a closed
graded trace on $\Omega^*(T,\CC\Gamma)$, then
$\eta$ defines
a cyclic cocycle on $C^\infty (T)\rtimes \Gamma$: we {\it assume }
that this cyclic cocycle {\it extends} to
a cyclic cocycle $\eta_{\mathcal{A}}$ on $\mathcal{A}$, then
$$
\langle {\rm ch}\, \Ind \D^{\rm sign} , \eta_{\mathcal{A}} \rangle = \int_{\widehat{X}/\Gamma}
L(T {\mathcal{F}} ) \wedge r^*(\Phi_\eta).
$$
where the index class is the one we defined in
$K_* (C(T)\rtimes_r \Gamma )$ and where $\Phi_\eta$
is the cohomology class in $H^* ((E\Gamma\times T)/\Gamma,\CC)$
associated to $\eta$ under the
Gorokhovsky-Lott isomorphism.

\begin{definition}
We shall say that a closed graded trace on $\Omega^*(T,\CC\Gamma)$,
is {\it holomophically extendable} if there exists a
dense holomophically closed subalgebra
$\mathcal{A}\subset C(T)\rtimes_r \Gamma$
with
$$C^\infty_c (T)\rtimes \Gamma
\subset
C^\infty (T,\B^\omega)\subset \mathcal{A}$$
and with the property that the cyclic cocycle defined by $\eta$ 
is cohomologous to a cocycle that extends
from $C^\infty_c (T)\rtimes \Gamma$ to $\mathcal{A}$.
\end{definition}

Making use of Proposition \ref{Prop}
in the previous section we thus obtain the following general

\begin{theorem} Let $\widehat{X}_\phi$ and $\widehat{X}_\psi$
two cut-and-paste equivalent $T\rtimes \Gamma$-proper manifolds
satisfying Assumption {\rm \ref{Lottgroupoid}}.
Let  $[c]\in
H^{*}( (E\Gamma \times T)/\Gamma ; \CC)$ be equal
to $\Phi_\eta$,
with $\eta$ a closed graded trace on $\Omega^*(T,\CC\Gamma)$
concentrated at the identity element.
If $\eta$ is holomorphically extendable, then
$$\int_{\widehat{X}_\phi/\Gamma}
L(T {\mathcal{F}}_\phi) \wedge r_\phi^*([c])\,=\,
\int_{\widehat{X}_\psi/\Gamma}
L(T {\mathcal{F}}_\psi ) \wedge s_\psi^*([c]).
$$
\end{theorem}

We shall now give examples where the assumptions
of the theorem are satisfied.

\subsection{Isometric actions.}

Assume that $\Gamma$ is Gromov Hyperbolic and preserves a Riemann
metric of $T$. Let $\omega \in \Omega^{{\rm dim}\, T -k}(T)$ be a
$\Gamma-$invariant differential form, then one defines a cyclic
cocycle $\tau_\omega$ on $C^\infty(T)$ by the formula: $$
\tau_\omega(f_0,f_1,\ldots , f_k) = \int_T f_0 d f_1 \wedge \ldots
\wedge d f_k \wedge \omega. $$ Let $\rho^\prime \in H^l(\CC
\Gamma, \CC)$ be a group cocycle and denote by $\tau_{\rho^\prime}
\in HC^l(\CC \Gamma)$ the associated cyclic cocycle constructed by
Connes (\cite{Co}). Consider then the cyclic cocycle $\tau^\prime$
of the algebraic tensor product $C^\infty(T) \rtimes \Gamma =
C^\infty_c(T \rtimes \Gamma)$ defined by: $$ \tau^\prime ( f_0
\gamma_0, f_1 \gamma_1, \ldots f_{k+l} \gamma_{k+l}) =
\tau_{\rho^\prime} \# \tau_\omega ( f_0 \otimes \gamma_0, \ldots ,
f_{k+l} \otimes \gamma_{k+l}) $$ where the $f_j \in C^\infty(T)$,
the $\gamma_j \in \Gamma$ and the $\#$ is defined in \cite{Co}
page 191. Then Jiang has proven 
the following (see Section 3 of \cite{Jiang}):
there exists a group cocycle $\rho$ cohomologous to
$\rho^\prime$  such that  $\rho$ has polynomial growth and $\tau =
\tau_\rho \# \tau_\omega$ extends as a cyclic cocycle to a dense
holomorphically closed
subalgebra $\mathcal{A}$ of $\Cr$. 
We conclude that $\tau$ is holomorphically extendable.
Let $\Phi_\rho \in
H^\ast ( ( E \Gamma \times T)/\Gamma ; \CC) $ denote the
associated cohomological class (see \cite{Go-Lo}) . Then
the Baum-Connes higher signature 
$$\int_{\widehat{X}/\Gamma} L(T {\mathcal{F}} )\wedge r^*(\Phi_\rho)$$
is a cut-and-paste invariant under our basic assumption \ref{Lottgroupoid}.
In fact, one can prove that
$r^*(\Phi_\rho)= f^*(\rho) \wedge [\pi^*(\omega)]$ where $f:
\widehat{X}/\Gamma \rightarrow B\Gamma$ denotes the classifying
map of the $\Gamma-$covering $\widehat{X} \rightarrow
\widehat{X}/\Gamma $ and $[\pi^*(\omega)]$ denotes the
differential form on $\widehat{X}/\Gamma$ induced by the
$\Gamma-$invariant differential form $\pi^*(\omega);$
thus, equivalently, we have proved that 
$\int_{\widehat{X}/\Gamma} L(T {\mathcal{F}} ) \wedge f^*(\rho)
\wedge [\pi^*(\omega)]$
is a cut-and-paste invariant if Assumption 
\ref{Lottgroupoid} is satisfied.

We should mention here that this example is automatically covered
by the results of subsection \ref{subsect:BC}
(assuming, as we do there, that the vertical tangent bundle
admits a $\Gamma$-equivariant spin structure). The next example,
on the other hand, is somewhat universal and it is {\it not}
covered by the Baum-Connes approach explained in subsection
\ref{subsect:BC}.

\subsection{The Godbillon-Vey signature $\sigma_{GV}$.}

Assume  $T=S^1$ and consider a $\Gamma-$equivariant fibration
$\pi_{\widehat{X}}: \widehat{X} \rightarrow S^1$ as in Subsection
\ref{subsect:Cut-and-paste}.

Let $X:=\widehat{X}/\Gamma$ and let $\mathcal{F}$
the induced foliation.
There exists a well-defined Godbillon-Vey
class $GV\in H^{3} (X,\RR)$.
The Godbillon-Vey signature is, by
definition $$\sigma_{GV}(X,\mathcal{F}):=\int_X L(T\mathcal{F})\wedge GV$$
By results of Connes \cite{Co} (see also the work of Moriyoshi-Nastume
\cite{Mor}),  we
know that there exists a dense and holomorphically closed
subalgebra $\mathcal{A}$
of $C(S^1)\rtimes_r \Gamma$ and a  closed graded trace
$\eta$ on $\Omega^* (S^1,\CC\Gamma)$ with the following properties:

\begin{itemize}
\item the cyclic cocycle defined by $\eta$ {\it extends}
to a cyclic cocycle $\eta_{\mathcal{A}}$ on $\mathcal{A}$
\item if $\Phi_{\eta}\in H^*
((E\Gamma\times S^1)/\Gamma,\RR)$ is the class corresponding to
$\eta$ under the Gorokhovsky-Lott isomorphism between the
homology of closed graded tra/-ces on $\Omega^*(S^1,\CC\Gamma)$
(concentrated at the identity)
and the cohomology of $(E\Gamma\times S^1)/\Gamma$, then
$GV=r^* \Phi_{\eta}$.
\end{itemize}

In other words, the so called Godbillon-Vey cyclic cocycle
is holomorphically  extendable and has the
property that $r^* \Phi_\eta$ is  equal
to the Godbillon-Vey class in $X$.

By applying the index theorem of Connes, we get
$$\langle\ch (\Ind (\mathcal{D}^{{\rm sign}})),\eta_{\mathcal{A}}\rangle=
\int_X L(T\mathcal{F})\wedge GV:=\sigma_{GV} (X,\mathcal{F}).
$$

It is important to remark that this formula holds with no assumption
on the group $\Gamma$ and its action on $T$.

\subsection{On the cut-and-paste invariance of $\sigma_{GV}$.}

Suppose now that $\widehat{X}_\phi$ and $\widehat{X}_\psi$
are two cut-and-paste equivalent $T\rtimes \Gamma$-proper manifolds
as in Subsection \ref{subsect:Cut-and-paste}.
By applying the above formula and Proposition
\ref{Prop},  we discover that
the Godbillon-Vey signature is a cut-and-paste invariant
if the middle-degree assumption \ref{Lottgroupoid} is satisfied.
More precisely:

\begin{theorem} Let $\widehat{X}_\phi$ and $\widehat{X}_\psi$
two cut-and-paste equivalent $T\rtimes \Gamma$-proper manifolds
satisfying Assumption {\rm \ref{Lottgroupoid}}. Let
$(X_\phi,\mathcal{F}_\phi)$ 
and $(X_\psi,\mathcal{F}_\psi)$ be the associated foliated manifolds.
Then
$$\sigma_{GV} (X_\phi,\mathcal{F}_\phi)= \sigma_{GV}
(X_\psi,\mathcal{F}_\psi)\,.$$
\end{theorem}

\subsection{Cut-and-paste invariance: 
the Baum-Connes approach.}\label{subsect:BC}

Now consider (see Connes \cite{Co}, page 114) the Baum-Connes rational
assembly map:
 $$
 \mu_\QQ: K_{*,\tau} ((E\Gamma\times T)/\Gamma)\otimes_\ZZ \QQ \rightarrow K_* ( \Cr)\otimes_\ZZ \QQ.
 $$
Here $\tau:= (E\Gamma\times T(T))/\Gamma$, with $T(T)$ denoting the tangent bundle
to $T$, and $K_{*,\tau} ((E\Gamma\times T)/\Gamma):= K_0 (B\tau,S\tau)$,
with $B\tau$ and $S\tau$ denoting the ball and sphere bundles of $\tau$.
In the general foliated case (\cite{Co}, Ch 2, Section 8.$\gamma$) there is a similar map
$$
 \mu_\QQ: K_{*,\tau} (BG)\otimes_\ZZ \QQ
\rightarrow K_0 ( C^*_r (M,\mathcal{F}))\otimes_\ZZ \QQ
 $$
with $BG$ the classifying space of the holonomy groupoid.
It is stated in work of Baum and Connes (\cite{Baum-Connes} page 12)
 that the rational injectivity of the latter assembly map implies
the leafwise homotopy invariance of the Baum-Connes higher signatures
of a foliation.
The line of reasoning is the following: if two foliations are
leafwise homotopy equivalent, then 
their leafwise
signature index classes are equal in
$K_0 (C^*_r (M,\mathcal{F}))$ (see \cite{Baum-Connes});
if moreover the Baum-Connes map is rationally injective,
then from the equality
of the index classes
 one deduces the equality of all the higher signatures.
The slogan here is the following: if the Baum-Connes map
is rationally injective , then {\it the equality of the index classes
implies the equality of all higher signatures. }

In the next theorem we shall prove 
the analogue of this fact in our groupoid $T\rtimes \Gamma-$setting.
In order to simplify our treatment we shall make the additional
assumption that the vertical tangent bundle $TZ$ to the fibration
$\widehat{X}\to T$, coming into the definition of $(T\rtimes \Gamma)$-proper
manifold, admits a $\Gamma$-invariant spin structure.

\begin{theorem}  Assume that the rational Baum-Connes map
$$
 \mu_\QQ: K_{0,\tau} ((E\Gamma\times T)/\Gamma)
\otimes_\ZZ \QQ \rightarrow K_0 ( \Cr)\otimes_\ZZ \QQ
 $$
is injective.
Let $\widehat{X}_\phi$ and $\widehat{X}_\psi$ be
two cut-and-paste equivalent $T\rtimes \Gamma$-proper manifolds
satisfying 
 Assumption {\rm \ref{Lottgroupoid}}
and such that the vertical tangent bundles both admit a $\Gamma$-invariant spin structure.
Then for any
$c \in  H^{\ast}( (E\Gamma \times T)/\Gamma ; \QQ)$ one has:
$$\int_{\widehat{X}_\phi
/{\Gamma} }
L( \widehat{X}_\phi
/{\Gamma})\wedge r_\phi^*(c)\,=\,\int_{\widehat{X}_\psi
/\Gamma }
L(\widehat{X}_\psi
/\Gamma )\wedge s_\psi^*(c).
$$
\end{theorem}
\begin{proof} 
We are going to show that for any
$c \in  H^{\ast}( (E\Gamma \times T)/\Gamma ; \CC)$ one has:
$$\int_{\widehat{X}_\phi
/{\Gamma} }
L( T \mathcal{F}_\phi )\wedge r_\phi^*(c)\,=\,\int_{\widehat{X}_\psi
/\Gamma }
L(T \mathcal{F}_\psi )\wedge s_\psi^*(c)
$$ which will prove the result. 
Let $$r_\phi: X_\phi\to (E\Gamma\times T)/\Gamma\;\;;\, s_\psi: X_\psi\to (E\Gamma\times T)/\Gamma$$
be the two classifying maps induced by 
$(\widetilde{\rho}_\phi,\pi_\phi):\widehat{X}_\phi\to E\Gamma\times T$
and $(\widetilde{\rho}_\psi,\pi_\psi):\widehat{X}_\psi\to E\Gamma\times T$ 
respectively (see subsection \ref{subsect:Cut-and-paste}).
Let $S_{\mathcal{F},_\phi}$, $S_{\mathcal{F},\psi}$ be the spinor bundles induced on the quotients
$X_\phi$, $X_\psi$ by the vertical spinor bundles.
Since the vertical tangent bundles carry a $\Gamma$-invariant Spin structure,
 we see that the bundles 
$$T X_\phi\oplus r_\phi^* \tau \;\;\;\text{and}\;\;\;T X_\psi\oplus s_\psi^* \tau$$
carry a Spin$^c$ structure.
 Then
$$[X_\phi, S^*_{\mathcal{F},_\phi}, r_\phi: X_\phi\to (E\Gamma\times T)/\Gamma]\,;\quad
[X_\psi, S^*_{\mathcal{F},_\psi}, s_\psi: X_\psi\to (E\Gamma\times T)/\Gamma]$$
define two geometric cycles in $K_{0,\tau} ((E\Gamma\times T)/\Gamma)$, see \cite{Co}
page 115. Notice, incidentally, that the Todd class is well
defined for a Spin$^c$ bundle (\cite{Co}, page 115). It follows from the very definition of the Baum-Connes map that
the image of $[X_\phi, S^*_{\mathcal{F},_\phi}, r_\phi]$
under $\mu$ is precisely equal to the index class $\Ind \D^{{\rm sign}}_\phi\in K_0 (\Cr)$ and 
similarly for 
$[X_\psi, S^*_{\mathcal{F},_\psi}, s_\psi]$.
From our Assumption \ref{Lottgroupoid} we know that the two index classes are equal; thus, by the assumed
rational injectivity of the Baum-Connes map,
we
have the
equality 
$$[X_\phi, S^*_{\mathcal{F},_\phi}, r_\phi]=
[X_\psi, S^*_{\mathcal{F},_\psi}, s_\psi]
\;\;\;\text{in}\;\;\;K_{0,\tau} ((E\Gamma\times T)/\Gamma)\otimes_\ZZ \QQ.$$
On the other hand, there is an isomorphism
$$K_{0,\tau} ((E\Gamma\times T)/\Gamma)\otimes_\ZZ \QQ\longrightarrow H_{{\rm even}} ((E\Gamma\times T)/\Gamma,\QQ)$$
which is given by the composition of the Chern character $\ch$ in $K$-homology
and the Thom isomorphism $\Phi$ in homology. Using Proposition 7 page 116 in \cite{Co}
we discover that 
$\Phi \circ \ch [X_\phi, S^*_{\mathcal{F},\phi}, r_\phi]$ equals
$$
(r_\phi)_*  \left( ( \ch (S^*_{\mathcal{F},\phi})\cdot
{\rm Td} (T X_\phi\oplus r^*_\phi \tau) \,)\cap [X_\phi] \right)\in H_{{\rm even}}
 ((E\Gamma\times T)/\Gamma,\QQ).
 $$
Using Lemma 4.4 of \cite{BGV} (pp. 148-150) we see that 
$$
(r_\phi)_*  \left( (\ch  S^*_{\mathcal{F},\phi} \cdot
{\rm Td} (T X_\phi\oplus r^*_\phi \tau  ) \,) \cap [X_\phi] \right)$$
equals
$$C \cdot
(r_\phi)_*  \left( ( L(T\mathcal{F}_\phi)\cdot
{\rm Td}_{\CC} (r^*_\phi \tau) \,)\cap [X_\phi] \right)$$
where $C \in\QQ^*$
 depends on the rank of $S^*_{\mathcal{F},\phi}$. Summarizing, it follows from our assumptions
that
$$\int_{X_\phi}
L( T \mathcal{F}_\phi )\wedge {\rm Td}_\CC (r^*_\phi \tau)\wedge r_\phi^*(c)\,=\,\int_{X_\psi}
L(T \mathcal{F}_\psi )\wedge {\rm Td}_\CC (s^*_\psi \tau) \wedge s_\psi^*(c).
$$
Applying this equality to $c:=({\rm Td}^{-1}_{\CC} \tau)\wedge b$
we get the equality
$$\int_{X_\phi}
L( T \mathcal{F}_\phi )\wedge r_\phi^*(b)\,=\,\int_{X_\psi}
L(T \mathcal{F}_\psi ) \wedge s_\psi^*(b).
$$
for each $b\in H^{*} ((E\Gamma\times T)/\Gamma,\QQ)$
which proves the theorem.

\end{proof}

\noindent {\bf Remark.} The Baum-Connes map $\mu_\QQ$ is injective many case,
we just mention three of them:

(a) $\Gamma$ has the Haagerup property (ie a-T-amenable) (see
Higson-Kasparov \cite{Hi-Ka}).

(b) $\Gamma$ is Gromov hyperbolic or more generally
$\Gamma$ is any discrete group acting properly by isometries on a weakly
bolic, weakly geodesic metric space of bounded coarse geometry
(see Kasparov-Skandalis \cite{Ka-Sk}).

(c) $\Gamma$ is a lattice in a semi-simple Lie group $G$
and $T=G/P$ where $P$ is a minimal parabolic subgroup of $G$
(for more results in this direction see Skandalis-Tu-Yu \cite{Sk-Tu-Yu}).

\newcommand{\foottext}[1]{{%
  \renewcommand{\thefootnote}{}\footnotetext{#1}}}

\foottext{Received by the editors January 31, 2004; revised May 11, 2004.}

\end{document}